\documentclass[12pt,english,fleqn]{article}
\usepackage[T1]{fontenc}
\usepackage[latin9]{inputenc}
\usepackage[a4paper]{geometry}
\usepackage{fancyhdr}
\usepackage{verbatim}
\usepackage{amsthm}
\usepackage{amsmath}
\usepackage{setspace}
\usepackage{amssymb}

\makeatletter

\providecommand{\tabularnewline}{\\}

\theoremstyle{plain}
\newtheorem{thm}{Theorem}
  \theoremstyle{remark}
  \newtheorem{rem}{Remark}
  \theoremstyle{plain}
  \newtheorem{lem}{Lemma}
  \theoremstyle{plain}
  \newtheorem{prop}{Proposition}
  \theoremstyle{plain}
  \newtheorem{cor}{Corollary}

\makeatother

\usepackage{babel}

\begin{document}
\newcommand{\myitemlabel}[1]{\mbox{\mdseries\upshape #1}} 
\newcommand{\ds}{\displaystyle } 

\global\long\def\myg{\gamma}
 \global\long\def\sP{\mathsf{P}}
 \global\long\def\sE{\mathsf{E}\,}
\global\long\def\RR{\mathbb{R}}

\global\long\def\MADR{\boldsymbol{R}}

\global\long\def\MADD{\boldsymbol{D}}

\global\long\def\md{d}

\global\long\def\mbd{\boldsymbol{d}_{N}}

\global\long\def\tn{t_{N}}

\global\long\def\myde{\delta_{N}}

\global\long\def\myal{\alpha_{N}}

\global\long\def\mybt{\beta_{N}}

\global\long\def\mLs{L_{N}}

\global\long\def\mM{M}

\global\long\def\Mp{S}
\global\long\def\rMp{\boldsymbol{\Mp}}

\global\long\def\mK{K_{N}}
\global\long\def\mko{k_{N}}
 \global\long\def\Id{\textrm{Id}}
 \global\long\def\myu{\boldsymbol{u}}
 \global\long\def\myb{\boldsymbol{b}}
 \global\long\def\OmKo{\left(\Id-\mK\right)^{-1}}
\global\long\def\Kmk{\left(\mK-\mko\right)}
 \global\long\def\Kmko{\Kmk^{-1}}
 \global\long\def\mqd{\textrm{q}_{2}}
 \global\long\def\mqo{\mbox{q}_{1}}
 \global\long\def\mqz{\textrm{q}_{0}}
 \global\long\def\onull{\left(\begin{array}{c}
 1\\
0\end{array}\right)}

\global\long\def\exO{e^{-\delta t\left(\Id-\mK\right)}}

\global\long\def\exSk{e^{-\delta t\left(1-\mko\right)}}

\global\long\def\PdcO#1{\frac{\Id+\delta t#1}{e^{\delta t}}}

\global\long\def\PdcSk#1{\frac{1+\delta t#1}{e^{\delta t}}}

\global\long\def\mywd{w_{2}}
 \global\long\def\myhd{h_{2}}

\global\long\def\lon{\lambda_{1,N}}

\global\long\def\mlon{\left|\lon\right|}

\global\long\def\ldn{\lambda_{2,N}}

\global\long\def\lin{\lambda_{i,N}}

\global\long\def\eon{\boldsymbol{e}_{1,N}}

\global\long\def\edn{\boldsymbol{e}_{2,N}}

\global\long\def\wdn{w_{2,N}}

\global\long\def\cy{\underline{y}}

\global\long\def\Rico{\MADR_{N}^{\,\mathrm{C}1}(\infty)}
\global\long\def\Dico{\MADD_{N}^{\,\!\mathrm{C}1}(\infty)}

\global\long\def\Ricd{\MADR_{N}^{\,\mathrm{C}2}(\infty)}
\global\long\def\Dicd{\MADD_{N}^{\,\mathrm{C}2}(\infty)}

\global\long\def\cdho{h_{R}}
 \global\long\def\cdhd{h_{D}}

\global\long\def\tm{\tau_{t}^{*}}

\global\long\def\JJ{J}
 \global\long\def\lo{\mbox{``\ensuremath{1}''}}
 \global\long\def\ld{\mbox{``\ensuremath{d}''}}
 \global\long\def\lV{\mbox{``\ensuremath{V}''}}
 \global\long\def\labt{\mbox{``\ensuremath{t}''}}
 \global\long\def\VG{\mathcal{V}_{G}}
 \global\long\def\EG{\mathcal{E}_{G}}
 \global\long\def\vD{v}

\global\long\def\matK{\mK}
 \global\long\def\ko{\mko}
 \global\long\def\qdv{\mqd}
 \global\long\def\qod{\mqo}
\global\long\def\oz{\onull}
 \global\long\def\zo{\left(\begin{array}{c}
 0\\
1\end{array}\right)}
 \global\long\def\bB{\myb}
\global\long\def\uU{\myu}

\global\long\def\emd#1{d^{(#1)}}
 \global\long\def\emf#1{f^{(#1)}}
 \global\long\def\emV#1{V^{(#1)}}
 \global\long\def\emR#1{R^{(#1)}}
 \global\long\def\emD#1{D^{(#1)}}

\global\long\def\dt{\Delta}
 \global\long\def\Ys#1#2{s_{#2}^{(#1)}}

\global\long\def\ao{a_{1}}
 \global\long\def\at{a_{2}}
 \global\long\def\aa{a}
 \global\long\def\gd{g_{2}}
 \global\long\def\go{g_{1}}
 \global\long\def\gG{g}
 \global\long\def\co{c_{1}}
 \global\long\def\cd{c_{2}}

\global\long\def\cR{C_{R}}
 \global\long\def\cD{C_{D}}
 \global\long\def\psR{\psi_{R}}
 \global\long\def\psD{\psi_{D}}
 \global\long\def\myfi{\phi}
 \global\long\def\mys{s}
 \global\long\def\mylR{l_{R}}
 \global\long\def\mylD{l_{D}}

\global\long\def\Fn#1#2{F_{N}(#1,#2)}
 \global\long\def\hH{H}

\global\long\def\Otn{O\left({\textstyle \frac{\tn}{N}}\right)}

\newcommand{\titleA}{Time Scales in Probabilistic Models\\ of Wireless Sensor Networks}

\title{\titleA  }

\date{February 28, 2013}

\author{Anatoly Manita%
\thanks{Department of Probability, Faculty of Mathematics and Mechanics, Lomonosov
Moscow State University, 119992, Moscow, Russia\protect \\
This is the extended presentation of results~\cite{Man-Msc-MM-2012}
announced at the GNEDENKO-100 Conference ({}``Probability Theory
and its Applications'', Moscow, June 26-30, 2012).\protect \\
The work is supported by Russian Foundation of Basic Research
(grant No.~12-01-00897).%
}}

\maketitle
~\vspace{-8.5ex}

\begin{center}
Lomonosov Moscow State University (Russia) \vspace{-0.95ex}

\par\end{center}

\begin{center}
\textsf{manita@mech.math.msu.su}
\par\end{center}
\begin{abstract}
We consider a stochastic model of clock synchronization in a wireless
network consisting of $N$ sensors interacting with one dedicated
accurate time server. For large $N$ we find an estimate of the final
time sychronization error for global and relative synchronization.
Main results concern a behavior of the network on different time scales
$\tn\rightarrow\infty$, $N\rightarrow\infty$. We discuss existence
of phase transitions and find exact time scales on which an effective
clock synchronization of the system takes place. 
\end{abstract}
\noindent \textbf{\small Keywords:}{\small{} Clock synchronization,
time to synchronization, timeliness in real-time systems, wireless
sensor networks, multi-dimensional Markov process, self-organization,
phase transitions}{\small \par}

\section{Introduction}

\label{sec:intro}

For many years distributed systems are a constant source of challenging
research problems. Clock synchronization is one of the most well-known
and widely discussed subjects. In distributed systems there is no
global clock. Clocks of different components tick at different rates.
But to achieve an efficient parallelization in their common job all
components of the distributed system need a common notion of time.
This is very important for real-time systems where scheduling and
timeliness play a crucial role (see, for example, \cite{FJNS-LSong-mk-2008,Zhu-Ma-Ryu-pat2009}).
Real-time applications must guarantee response within strict time
constraints. Similarly, Wireless Sensor Networks (WSNs), having their
own specific character among distributed systems, usually meet many
real-time requirements. In particular, energy saving algorithms in
WSNs would be impossible without a global consistent timescale~\cite{energy-conserv-Wu-2005}.
Nowadays complex real-time embedded systems are distributed and have
many active components. Modern WSNs also consist of very large number
of sensors (\cite{Song-Schott-Large-SN,Austr4aut2010}). Hence an
important requirement for clock synchronization in distributed systems
is the scalability of corresponding algorithms. 

In this paper we consider a mathematical model of a large WSN communicating
with one dedicated server of accurate time. The model is based on
a multi-dimensional stochastic process. The network consists of $N$
sensors equipped with non-perfect clocks. By using message timestamps
sensors share information about their local times and receive accurate
time from the dedicated server. It is assumed that all messages are
sent at random time moments. Main focus of the paper is the asymptotical
behavior of the network in the limit when the number~$N$ of sensors
tends to infinity.

This article is organized as follows. The model is defined in Section~\ref{sec:st-model-Ns}.
It depends on a few parameters. The network is assumed to be symmetric.
Here we do not touch many problems (synchronization protocols, energy
saving optimization etc.) which are very important for practical use
of WSNs. Such questions were widely discussed by many authors \cite{SundBuyKshem-review,Li-Rus,energy-conserv-Wu-2005,Song-Schott-Large-SN,Romer}.
The present article is a mathematical paper. We intentionally make
our model as transparent as possible to keep this paper short and
to illuminate main results about existence of different phases in
the evolution of the network (Theorems~1-3 in Section~\ref{sec:main-res-WSN-serv}).
It appears that the proposed model has almost explicit solution: in
all asymptotical expressions we can play with parameters and discover
many interesting phenomena related with large networks. In Section~\ref{sec:future-w}
we discuss possible generalization of the present results to more
general models. We do not seek for direct practical implementation
of our results. But we believe that such results are interesting not
only from theoretical viewpoint. 

Our probablistic technique is similar to that was recently used to
study mathematical models of multi-processor parallel computing~\cite{MitraMitrani,Bert-Tsits,Francois-Arx,Francois-LNCS}
and collective behavior in abstract particle systems with synchronization-like
interaction~\cite{man-Shch,mal-man-TVP,manita-umn,Malyshkin,manita-TVP-large-ident-eng,Man-F-T,Man-Rome}.
Our stochastic models can be considered also as special classes of
self-organizing systems~\cite{Bettstetter-SO-CN,Dressler-Actor,Pikov-Synchro,Strogatz-SYNC}.

\section{Mathematical model}

\label{sec:st-model-Ns}

There is a wireless network with $N+1$ nodes. Each node~$i$ has
its own clock. Let $x_{i}(t)$ be the value of clock~$i$ (the \emph{local
time} at~node~$i$). The physical time is denoted by $t\in\mathbb{R}_{+}$.
Nodes $2,3,\ldots,N+1$ correspond to \emph{sensors}. We assume that
clocks of the sensors are identical but \emph{not perfect} and progress
in the following way:\begin{equation}
x_{j}(t)=x_{j}(0)+vt+\sigma B_{j}(t),\quad j=2,\ldots,N+1\,.\label{MAD:eq:x-vt-sB}\end{equation}
Node~$1$ is a \emph{time server}, we assume that its clock is \emph{perfect.}
Clock~1 reports a local time $x_{1}(t)$ which is a linear function
of the physical time:\begin{equation}
x_{1}(t)=x_{1}(0)+rt\,.\label{MAD:eq:x-rt}\end{equation}
The constants $r$ and $v$ are positive, they are called \emph{fre\-quen\-cies}
of the corresponding clocks~(see, for example, \cite{SundBuyKshem-review}).
In general, these frequencies are not equal: $r\not=v$. The parameter~$\sigma>0$
in~(\ref{MAD:eq:x-vt-sB}) corresponds to the strength of a random
noise related with imperfect clock of a sensor. For simplicity the
random noise in~(\ref{MAD:eq:x-vt-sB}) is modelled with independent
standard Brownian motions \[
\left(B_{j}(t),\,\, t\geq0\right),\qquad j=2,\ldots,N+1.\]
This white noise assumption is usual for many oscillator clock models~\cite{Zhu-Ma-Ryu-pat2009,Austr4aut2010}.

As it was mentioned in Section~\ref{sec:intro} (see also~\cite{SimSpagBarStrog})
the problem of clock synchronization is critical for the proper work
of WSNs. To synchronize local times of a pair of sensors the most
natural solution is to send time-stamped messages between them. After
reading a new coming message the receiver ajusts its clock to the
local time of the sender recorded in this message. Once being synchronized
the clocks of the pair will diverge immediately due to assumptions
made in~(\ref{MAD:eq:x-vt-sB}) and~(\ref{MAD:eq:x-rt}). 

We continue with the formal definition of the model. We shall call~(\ref{MAD:eq:x-vt-sB})--(\ref{MAD:eq:x-rt})
a \emph{free dynamics}. The free dynamics generates independent evolutions
at nodes of the network. Now we add some special \emph{interaction}
between nodes. Namely, at random time moments each node sends messages
to other nodes. Our assumptions are
\begin{itemize}
\item with the rate $\alpha>0$ the \emph{server node}~1 generates a message
containing information about the current value of $x_{1}$ and sends
this message to one of the \emph{sensors} which is chosen randomly
with probability~${\displaystyle \frac{1}{N}}$~;
\item independently \emph{each sensor}, with the rate $\beta>0$, generates
a message about its local time and sends it to another \emph{sensor}
which is chosen randomly with probability~${\displaystyle \frac{1}{N-1}}$~;
\item messages reach their destinations instantly (there are no transmission
delays);
\item if sensor~$j$ receives a~mes\-sage from sensor~$i$ at some (random)
time $\tau$ then clock~$j$ is immediately ajusted to the value
of clock~$i$:~ $x_{j}(\tau+0)=x_{i}(\tau)$;
\item between receiving of subsequent messages sensor nodes evolve according
to the free dynamics~(\ref{MAD:eq:x-vt-sB}) ;
\item server node~1 always follows the free dynamics~(\ref{MAD:eq:x-rt}).
\end{itemize}
Let us remind some terminology used above: a sequence of events $0=\tau_{0}<\cdots<\tau_{n}<\cdots$
generated with rate~$\delta>0$ is called also a Poisson flow of
intensity~$\delta$. It means that intervals between events $\tau_{n+1}-\tau_{n}$
are independent exponentially distributed random variables with mean~$\myg^{-1}$:
\begin{equation}
\sP\left(\tau_{n+1}-\tau_{n}>s\right)=e^{-\delta s},\quad\quad\sE\left(\tau_{n+1}-\tau_{n}\right)=\delta^{-1}\,.\label{eq:exp-distrub}\end{equation}

Hence we defined the multi-dimensional stochastic process \[
x(t)=\left(x_{1}(t),x_{2}(t),\ldots,x_{N+1}(t)\right),\quad\quad t\geq0,\]
which appears to be a continuous time Markov process with values in~$\mathbb{R}^{N+1}$.
This process is an idealized mathematical model of the homogeneous
wireless sensor network with dedicated accurate time server.

It should be noted that non-Markovian models of WSNs can be also considered
in the framework of the present paper. We postpone discussion on possible
generalizations to Subsection~\ref{sub:many-notation}.

To estimate desynchronization in the network it is convenient to consider
the following functions on the configuration space~$\mathbb{R}^{N+1}$:
\begin{eqnarray}
R(x) & := & \frac{1}{N}\,\sum_{j=2}^{N+1}\left(x_{j}-x_{1}\right)^{2},\label{eq:Rx}\\
D(x) & := & \frac{1}{(N-1)N}\,\sum_{2\leq j_{1}<j_{2}}\left(x_{j_{2}}-x_{j_{1}}\right)^{2}.\label{eq:Dx}\end{eqnarray}
The first function corresponds to deviations of sensor's clocks from
the accurate time~$x_{1}$. The second function is related only to
\emph{internal} inconsistency between the sensor's local times $x_{2},\ldots,x_{N+1}$.
This function $D(x)$ may be more useful for causuality-based real-time
models where the right ordering of events is more significant than
the global synchronization to physical time (see, for example, \cite{zhao-liu-lee-2007}).
Since the both functions~$R(x)$ and $D(x)$ are averagings of $(x_{k}-x_{j})^{2}$,
the \emph{squares} of offsets between clocks, the true sense of a
time synchronization error have their square roots $\sqrt{R(x)}$
and $\sqrt{D(x)}$. The physical dimension of $\sqrt{R(x)}$ and $\sqrt{D(x)}$
is \emph{time}, i.e., they are measured in \emph{seconds}. 

Note that for any $t$ values $R(x(t))$ and $D(x(t))$ are random.
It is more convenient to deal with their \emph{expectations}:

\begin{center}
$\MADR_{N}(t):=\mathsf{E}\, R(x(t))\,\,$ and $\,\,\MADD_{N}(t):=\mathsf{E}\, D(x(t))\,.$ 
\par\end{center}

\section{Main results}

\label{sec:main-res-WSN-serv} Let $x(t)$ be the Markov process introduced
in Section~\ref{sec:st-model-Ns}.

In Subsections~\ref{sub:N-fixed-t-infty}--\ref{sub:big-initial}
we present various results on asymptotical behavior of the functions
$\MADR_{N}(t)$ and $\MADD_{N}(t)$. These results will provide us
with detailed information on collective behavior of the network in
the limit when $N\rightarrow\infty$. Here we give \emph{one} \emph{particular
corollary} of our theorems which is a good \emph{illustration} to
an approach followed in this paper.\smallskip{}

\noindent \textbf{Corollary.} Assume that $v\not=r$. For simplicity
take $x_{1}(0)=\cdots=x_{N}(0)=0$ as initial state of the network.
Let $\mys>0$ and $\gamma>0$ be parameters of a new time scale. Namely,
we put $t=\mys N^{\gamma}$ and look after the function $\MADD_{N}(t)$
when $N\rightarrow\infty$. It appears that \[
\MADD_{N}(\mys N^{\gamma})\sim C(\mys,\gamma)\, N^{\phi(\gamma)}\]
where

\[
\phi(\gamma)=\begin{cases}
\gamma, & \gamma\leq\frac{1}{2},\\
3\gamma-1, & \frac{1}{2}<\gamma\leq1,\\
2, & \gamma>1.\end{cases}\]
The function $C(\mys,\gamma)=C(\mys,\gamma;\sigma,\alpha,\beta)>0$
depends on parameters $\sigma,\alpha,\beta$ and is increasing in~$\mys$
for fixed~$\gamma$. \smallskip{}

The choice $t=\mys N^{\gamma}$ means that we consider a new time
unit $N^{\gamma}$ and the {}``size'' of configuration $x_{2},\ldots,x_{N}$
on this scale is of order~$N^{\phi(\gamma)/2}$. The conclusion is
that on different time scales the large network shows very different
types of behavior. Note that the function~$\phi(\gamma)$ is not
smooth. This situation is very like the \emph{phase transitions} phenomena
in models of statistical physics. We give more details on time scales
in~Subsection~\ref{sub:time-sc-analysis}. \smallskip{}

To understand in what sense the algorithm of Section~\ref{sec:st-model-Ns}
drives the network to synchronization one should answer several questions.

\subsection{$N$ is fixed, $t\rightarrow\infty$}

\label{sub:N-fixed-t-infty} The first question is a long-time behavior
of the stochastic process $$x(t)=\left(x_{1}(t),\ldots,x_{N+1}\right).$$
Since $r>0$ and $v>0$ it is quite clear from definition of the model
that $x_{j}(t)\rightarrow\infty$ ($t\rightarrow\infty$) for each
node $j$. To have more detailed information on the process we put
a moving observer at the point $x_{1}(t)$. From the viewpoint of
this observer states of the sensor nodes $2$, $\ldots$, $N+1$ are
given by the vector \begin{equation}
\cy(t)=\left(y_{2}(t),\ldots,y_{N+1}(t)\right),\quad y_{j}=x_{j}-x_{1}.\label{eq:yj-xj-x1}\end{equation}

It follows from the general theory of Markov processes (see Doeblin
condition in~\cite{Doob53}) that $\cy(t)$ has a limiting distribution
as $t\rightarrow\infty$. Similar situations were discussed in~\cite{manita-umn,Man-F-T}
for different synchronization models. Hence our\emph{ first conclusion}
is: ~ for fixed parameters $\sigma$, $\alpha>0$, $\beta>0$ and~$N$
the Markov process $\cy(t)$ is ergodic and converges to its equlibrium
as $t\rightarrow\infty$. This limiting distribution (which is a probability
measure on $\mathbb{R}^{N}$) can not be obtained in an explicit form.
Nevertheless, some important functionals of the limiting distribution
can be found explicitly. Note that $R(x)=R(y)$ and $D(x)=D(y)$ where
$x$ and $\cy$ are related by~(\ref{eq:yj-xj-x1}). Thus the clocks~$x_{i}(t)$
of sensor nodes are synchronized in the following sense: mean values
of time synchronization errors stabilize and do not change in time
anymore. Namely, assuming that $N$ is fixed and $t\rightarrow\infty$
we can prove the following statement.
\begin{thm}
\label{thm:N-fixed-t-to-inf}~We distinguish two cases: the \textbf{skew}
$v-r$ of a sensor's clock relative to the accurate time is \textbf{zero}
or \textbf{nonzero}.
\begin{description}
\item [{Case~1:}] Assume that $v=r$. Then \begin{eqnarray*}
\MADR_{N}(t) & \rightarrow & \Rico:=\frac{\sigma^{2}}{\alpha}\, N\,,\\
\MADD_{N}(t) & \rightarrow & \Dico\,\,\sim\,\,\frac{2\,\sigma^{2}}{\alpha+\beta}\, N\,\,.\end{eqnarray*}

\item [{Case~2:}] Assume that $v\not=r$. Then \begin{eqnarray*}
\MADR_{N}(t) & \rightarrow & \Ricd:=\frac{\,\left(v-r\right)^{2}}{\alpha^{2}}\, N^{2}\,.\\
\MADD_{N}(t) & \rightarrow & \Dicd\sim\frac{\,\left(v-r\right)^{2}}{\alpha\left(\alpha+\beta\right)}\, N^{2}\,.\end{eqnarray*}

\end{description}
\end{thm}
As usual we write here and below \qquad\(
{\displaystyle f_{N}\sim g_{N}\,\quad\mbox{iff}\quad\,\lim_{N}\,\frac{f_{N}}{g_{N}}=1}.\)

Now we come to a natural and \emph{very important question}: how many
time we should wait until the networks of {}``size''~$N$ will
be synchronized~? We devote to this question the next Subsections~\ref{sub:phases-1-3}--\ref{sub:res-collect-displacement}.

Comparing Case~1 and Case~2 in the above theorem we see that result
of the global synchronization of sensor's clocks is worse in the biased
case $v\not=r$. Indeed,  condition $v\not=r$ means that there is
a systematic error in client's clocks which appears to be more essential
(a resulting {}``clock offset'' is of order $\sqrt{\Ricd}=O(N)$~)
than in situation of a pure random noise errors under condition~$v=r$
(when the {}``clock offset'' is of order $\sqrt{\Rico}=O(N^{1/2})$~).

Note also that in Case~2 the limiting values $\Ricd$ and $\Dicd$
do not depend on~$\sigma$, the noise parameter of the sensor's clocks.
Moreover, these values give the right asymptotics even for $\sigma=0$.
This means that the so called locked synchronization~\cite{SimSpagBarStrog}
is not possible in our model due to the stochastic nature of the message-passing
algorithm.

\subsection{$N\rightarrow\infty$, $t\rightarrow\infty$. Phases of synchronization}

\label{sub:phases-1-3} Synchronization in networks with large (or
growing) number of nodes is of special interest. Our approach is to
consider limits when both the number of client nodes~$N$ and the
time $t$ grow to infinity. More precisely, we shall consider sequences
$(N,\tn)$, where physical time $t=\tn$ is some increasing function
of~$N$. We shall see below that for different choices of $\tn$
network exhibits different asymptotical behavior. One can say that
a large network passes different phases on its road to synchronization.
For our model there exist \emph{several time scales} $t=\tn$ of qualitatively
different behavior. We always assume that $N\rightarrow\infty$. 

To begin \emph{we assume} in the present subsection that initial distribution
of the local clocks is such that sequences $\left\{ \MADR_{N}(0)\right\} $
and $\left\{ \MADD_{N}(0)\right\} $ are bounded in~$N$: \begin{equation}
\sup_{N}\MADR_{N}(0)<+\infty,\qquad\sup_{N}\MADD_{N}(0)<+\infty.\label{eq:R0D0}\end{equation}
In Subsec.~\ref{sub:big-initial} we shall discuss the situation
when $\left\{ \MADR_{N}(0)\right\} $ and $\left\{ \MADD_{N}(0)\right\} $
are unbounded. As in Subsect.~\ref{sub:N-fixed-t-infty} we present
results separately for Case~1 and Case~2.
\begin{thm}
\label{thm:v-eq-r-3-phases}Let assumption~(\ref{eq:R0D0}) hold.
\begin{description}
\item [{Case~1:}] $v=r$ --- zero skew. 
\end{description}
\renewcommand{\theenumi}{\mbox{\upshape P\arabic{enumi}}}
\begin{enumerate}
\item ~ $\tn/N\rightarrow0$ ~ (phase of initial desynchronization):\[
\left(\begin{array}{c}
\MADR_{N}(\tn)\\
\MADD_{N}(\tn)\end{array}\right)\sim\left(\begin{array}{c}
\sigma^{2}\\
2\sigma^{2}\end{array}\right)\,\tn\]

\item ~ $\tn/N\rightarrow c>0$ ~ (phase of effective syn\-chro\-ni\-zation):\[
\left(\begin{array}{c}
\MADR_{N}(\tn)\\
\MADD_{N}(\tn)\end{array}\right)\sim\,(-cM)^{-1}\left(\Id-e^{cM}\right)\,\left(\begin{array}{c}
\sigma^{2}\\
2\sigma^{2}\end{array}\right)\,\tn\,,\]
 where $\Id$ is the identity map in~$\RR^{2}$, \begin{equation}
M:=\left(\begin{array}{cc}
-\alpha & \,\,0\\
2\alpha & \,\,-2(\alpha+\beta)\end{array}\right)\,.\label{eq:matr-M}\end{equation}

\item ~ $\tn/N\rightarrow+\infty$ ~ (phase of final stabilization): \[
\left(\begin{array}{c}
\MADR_{N}(\tn)\\
\MADD_{N}(\tn)\end{array}\right)\sim\,\left(\begin{array}{c}
\sigma^{2}/\alpha\\
2\sigma^{2}/(\alpha+\beta)\end{array}\right)\, N\,\sim\left(\begin{array}{c}
\Rico\\
\Dico\end{array}\right).\]

\end{enumerate}
\end{thm}
\begin{rem}
By using definition~(\ref{eq:matr-M}) of the matrix~$M$ it is
easy to check that $\mM$ has two distinct negative eigenvalues: $\lambda_{1}=-\alpha,$
$\lambda_{2}=-2(\alpha+\beta)$. Moreover, $(-cM)^{-1}\left(\Id-e^{cM}\right)\rightarrow\Id$
as $c\rightarrow+0$, $e^{cM}\rightarrow0$ as $c\rightarrow+\infty$
and \[
\,(-M)^{-1}\,\left(\begin{array}{c}
1\\
2\end{array}\right)=\left(\begin{array}{c}
1/\alpha\\
2/(\alpha+\beta)\end{array}\right).\]
 
\end{rem}
From this remark we see that phase~P2 joins smoothly asymptotics
of phases~P1 and~P3.
\begin{rem}
\label{rem:v=00003Dr-stages}One can easily calculate $\MADR_{N}(t)$
and $\MADD_{N}(t)$ in the simplest case $v=r$, $\alpha=\beta=0$
(no bias, no synchronizing interaction). It appears that $\MADR_{N}(t)=\sigma^{2}t$
and $\MADD_{N}(t)=2\sigma^{2}t$. This observation brings us to the
following explanation of the phase~1: on time intervals $[0,o(N)]$
the influence of synchronizing jumps (clock adjustments) is negligible
with respect to the impact of the random noise of the free dynamics.
On the time scale of phase~P2 $(\tn=cN$) the cumulative effect of
individual synchronizing adjustments become of the same order as the
effect of random noise, so we obtain the non-trivial dependence of
$\MADR_{N}(\tn)$ and $\MADD_{N}(\tn)$ on~$c$. On~phase~P3 the~values
of $\MADR_{N}(\tn)$ and $\MADD_{N}(\tn)$ correspond to the synchronized
network. It means that the \emph{effective synchronization} takes
place on the times of order~$N$.\end{rem}
\begin{thm}
\label{thm:v-NOT-eq-r-3-phases} Assume that condition~(\ref{eq:R0D0})
holds.
\begin{description}
\item [{Case~2:}] $v\not=r$ --- nonzero skew.
\end{description}
\renewcommand{\theenumi}{\mbox{\upshape P\arabic{enumi}}}
\begin{enumerate}
\item ~ $\tn/N\rightarrow0$ ~ (phase of initial desynchronization):\[
\MADR_{N}(\tn)\sim\frac{1}{2}\,(v-r)^{2}\,\tn^{2}\,.\]
\begin{description}  \item[\myitemlabel{P1a.  }]  If ~~$\ds\frac{\tn}{\sqrt{N}}\,\rightarrow0$,
then $\MADD_{N}(\tn)\sim2\sigma^{2}\tn$. \item[\myitemlabel{P1b.  }] If\textup{
~~$\ds\frac{\tn}{\sqrt{N}}\,\rightarrow\co,\quad\co>0,$ }~ then\begin{eqnarray*}
\MADD_{N}(\tn) & \sim & \left(2\sigma^{2}+\frac{1}{3}\alpha(v-r)^{2}\,\co^{2}\right)\tn\,.\end{eqnarray*}
\item[\myitemlabel{P1c.  }] If ~~$\ds\frac{\tn}{\sqrt{N}}\,\rightarrow\infty$~
~but~ ~$\ds\frac{\tn}{N}\,\rightarrow0$,~ then~ $\MADD_{N}(\tn)\sim\frac{1}{3}\alpha(v-r)^{2}\,\tn^{3}/N$.
\end{description} 
\item ~ $\tn/N\rightarrow c>0$ (phase of effective syn\-chro\-ni\-zation):
there exist functions $\cdho$ and $\cdhd$ not depending on~$N$
such that\[
\left(\begin{array}{c}
\MADR_{N}(\tn)\\
\MADD_{N}(\tn)\end{array}\right)\sim\,\left(\begin{array}{c}
\cdho(c)\\
\cdhd(c)\end{array}\right)\,\left(v-r\right)^{2}\,\tn^{2}\]

\item ~ $\tn/N\rightarrow+\infty$ ~ (phase of final stabilization): \[
\left(\begin{array}{c}
\MADR_{N}(\tn)\\
\MADD_{N}(\tn)\end{array}\right)\sim\left(\begin{array}{c}
{\displaystyle \frac{\,\left(v-r\right)^{2}}{\alpha^{2}}}\\
\,\\
{\displaystyle \frac{\,\left(v-r\right)^{2}}{\alpha\left(\alpha+\beta\right)}}\end{array}\right)\, N^{2}\,\sim\left(\begin{array}{c}
\Ricd\\
\,\\
\Dicd\end{array}\right).\]

\end{enumerate}
\end{thm}
\begin{rem}
Explicit form of the function $\cdho$ and $\cdhd$ is given in~(\ref{eq:h1-h2-expl-form}),
see Section~\ref{sec-proofs-Ths}. It appears that \[
\left(\begin{array}{c}
\cdho(c)\\
\cdhd(c)\end{array}\right)\sim\left(\begin{array}{c}
1/2\\
\alpha c/3\end{array}\right)\quad(c\rightarrow+0),\]
\[
\left(\begin{array}{c}
\cdho(c)\, c^{2}\\
\cdhd(c)\, c^{2}\end{array}\right)\,\rightarrow\,\left(\begin{array}{c}
\alpha^{-2}\\
\alpha^{-1}(\alpha+\beta)^{-1}\end{array}\right)\quad(c\rightarrow+\infty).\]
Hence the phase~P2 {}``continuously'' joins asymptotics of phases~P1c
and~P3.
\end{rem}

\subsection{Situation of a {}``big initial disorder''}

\label{sub:big-initial}

As in preceeding subsection we consider here different time scales
$t=\tn$ assuming that $\tn\rightarrow\infty$ as $N\rightarrow\infty$.
Our goal is to discuss what happens with sychronization phases if
the sequences $\MADR_{N}(0)$ and $\MADD_{N}(0)$ grow to infinity
as $N\rightarrow\infty$. This problem is not too complicated and
can be solved explicitely for any concrete assumption about initial
disorder, i.e., about increasing rate of $\MADR_{N}(0)$ and $\MADD_{N}(0)$.
In this paper we shall not give universal and general answer to this
question since it cannot be presented in short and transparent expressions.
To~avoid cumbersome formulae we just describe the answer in a few
general words. There exist a so called \emph{initial disorder decay}
interval $(0,\tn^{\circ})$. The function $\tn^{\circ}$ depends on
$\MADR_{N}(0)$, $\MADD_{N}(0)$ and all parameters of the model.
Under above assumptions $\tn^{\circ}\rightarrow\infty$ as $N\rightarrow\infty$.
If~$\tn=o(\tn^{\circ})$, then the values $\MADR_{N}(0)$ and $\MADD_{N}(0)$
enter in asymptotics of $\MADR_{N}(\tn)$ and $\MADD_{N}(\tn)$, but
they don't enter in these asymptotics for~$\tn$ such that $\tn/\tn^{\circ}\rightarrow\infty$. 

Given a function~$\tn^{\circ}$, one should compare orders of $\tn^{\circ}$,
$N^{1/2}$ and $N$ to see intersections of the \emph{initial disorder
decay phase} (IDDP) with the phases~P1, P2 and~P3 of Subsection~\ref{sub:phases-1-3}.
The IDDP dominates over any of the phases~P1 (P1a--P1c), P2 and~P3.
So~depending on the initial disorder in the network one can see the
following different sequences of consequtive phases in evolution of~the~network:
IDDP--P1--P2--P3, IDDP--P1b--P1c--P2--P3, IDDP--P1c--P2--P3, IDDP--P2--P3
or IDDP--P3.

\subsection{Collective displacement from the etalon time}

\label{sub:res-collect-displacement}

In Subsections~\ref{sub:N-fixed-t-infty}--\ref{sub:big-initial}
we studied $R(x)$ and $D(x)$ which are quadratic functions of the
clock configuration~$x\in\RR^{N+1}$. For completeness we consider
here a linear function $\md(x)$ which gives some information on a
collective displacement of the sensors clocks~$x_{2},\ldots,x_{N}$
from the clock of the server~$x_{1}$. Define functions

\begin{equation}
\md(x)=\frac{1}{N}\,\sum_{j=2}^{N+1}x_{j}-x_{1}\,,\qquad\md:\,\mathbb{R}^{N+1}\rightarrow\mathbb{R}^{1},\label{eq:d-x-def}\end{equation}
and \[
\mbd(t)=\mathsf{E}\,\md(x(t))\,,\qquad\mbd:\,\RR_{+}\rightarrow\RR.\]

\begin{thm} \label{th:dN-asympt}

Assume that $\sup_{N}\mbd(0)<+\infty$. 
\begin{description}
\item [{Case~1:}] $v=r$ --- zero skew. It appears that $\mbd(t)\rightarrow0$
as $t\rightarrow\infty$ for any fixed~$N$. 
\item [{Case~2:}] $v\not=r$ --- nonzero skew. For fixed~$N$ \[
\mbd(t)\rightarrow(v-r)\alpha^{-1}N\quad\mbox{as}\quad t\rightarrow\infty.\]
For time scales $\tn\rightarrow\infty$ ($N\rightarrow\infty$) the
following statements hold:

$\circ$~ if $\tn/N\rightarrow0$, then $\mbd(\tn)\sim(v-r)\tn$~,

$\circ$~ if $\tn/N\rightarrow c>0$, then \[
\mbd(\tn)\sim\left(1-\exp\left(-\alpha c\right)\right)\,(v-r)\alpha^{-1}N\,\sim\,\frac{1-\exp(-\alpha c)}{\alpha c}\cdot(v-r)\,\tn\,,\]

$\circ$~ if $\tn/N\rightarrow+\infty$, then $\mbd(\tn)\sim\,(v-r)\alpha^{-1}N\,$~.

\end{description}
\null \end{thm}

The proof of this theorem is similar to proofs of Theorems~\ref{thm:N-fixed-t-to-inf}--\ref{thm:v-NOT-eq-r-3-phases}
(see Section~\ref{sec:proofs}) but technically it is much easier.
So we omit it.
\begin{rem}
Let us consider the degenerated model with $\alpha=0$. It is easy
to prove that \begin{equation}
\mbd(\tn)-\mbd(0)=(v-r)\tn\,.\label{eq:d-t-alf-0}\end{equation}
Comparing this with Theorem~\ref{th:dN-asympt} we see that the formal
limit $\alpha\rightarrow0$ does not turn results of Theorem~\ref{th:dN-asympt}
into~(\ref{eq:d-t-alf-0}). Hence we conclude that limits $N\rightarrow\infty$
and $\alpha\rightarrow0$ do not commute.
\end{rem}

\section{Proofs}

\label{sec:proofs}

\subsection{Conditional averaging}

Let $0=\tau_{0}<\cdots<\tau_{n}<\cdots$ be sequence of time moments
when messages are sent (received). It follows from definition of the
model (Section~\ref{sec:st-model-Ns}) that $\left\{ \tau_{m}-\tau_{m-1}\right\} _{m=0}^{\infty}$
are independent i.d. r.v. having exponential distribution with mean
$\left(\alpha+N\beta\right)^{-1}$. In the sequel we will refer to
this condition on $\left\{ \tau_{n}\right\} $ as to a \textbf{\emph{Markovian
assumption}}. 

Let $\Pi_{t}=\max\left\{ m:\,\tau_{m}\leq t\right\} $.  To get $\MADR_{N}(t)$
and $\MADD_{N}(t)$ we will calculate the chain of conditional expectations
as follows \[
\sE\left(\cdot\right)=\sE\,\left(\sE\left(\sE\left(\cdot\,|\,\left\{ \tau_{j}\right\} _{j=1}^{\infty}\right)\,|\,\Pi_{t}\right)\right).\]

Let $f=f(x)$ be some function on the configuration space~$\mathbb{R}^{N+1}$.
Introduce notation \begin{equation}
\emf n=\sE\left(f(x(\tau_{n}))\,|\,\left\{ \tau_{j}\right\} _{j=1}^{\infty}\right),\qquad n=1,2,\ldots\,.\label{eq:fn-f-tau-j}\end{equation}
 Hence $\emf n$ is a random variable functionally depending on the
sequence $\left\{ \tau_{j}\right\} _{j=1}^{\infty}$. In other words,
to any function $f=f(x)$ we put in correspondence $f\mapsto\left\{ \emf n\right\} $
the sequence of random variables $\left\{ \emf n\right\} $. 

Denote \[
V(x)=\left(\begin{array}{c}
R(x)\\
D(x)\end{array}\right),\qquad V:\,\mathbb{R}^{N+1}\rightarrow\mathbb{R}^{2},\]
where functions $R(x)$ and $D(x)$ are defined in~(\ref{eq:Rx})--(\ref{eq:Dx}).
Recall definition of the function~$d(x)$ from~(\ref{eq:d-x-def}).
In the next subsection we study random sequences $\left\{ \emV n\right\} $
and $\left\{ \emd n\right\} $ corresponding to the above functions~$V(x)$
and $d(x)$ according to~(\ref{eq:fn-f-tau-j}).

\subsection{Recurrent equations}

\label{sub:many-notation}

In this subsection we obtain difference equations for $\emV n=\left(\begin{array}{c}
\emR n\\
\emD n\end{array}\right)$ and $\emd n$. To present this result we need some convenient notation
which will be used till the end of the paper. Namely, we introduce
the following reals\begin{equation}
\begin{array}{rclcrclcrclcc}
\myb & = & v-r, & \quad & \myal & = & {\displaystyle \frac{\alpha}{N}\,,} & \quad & \mybt & = & {\displaystyle \frac{\beta}{N-1}\,,} & \quad & \quad\\
\\\myu & = & 2\myb\,, & \quad & \myde & = & \alpha+N\beta, & \quad & \mko & = & 1-\myde^{-1}\myal\,,\end{array}\label{eq:notation-scalars}\end{equation}
and two-dimensional vectors \[
\mqo=\sigma^{2}\left(\begin{array}{c}
1\\
2\end{array}\right),\quad\mqd=(v-r)^{2}\onull,\quad\mqz=\myu\onull.\]
We consider also a $2\times2$-matrix \begin{equation}
\mLs=\left(\begin{array}{cc}
-\myal & 0\\
2\myal & -2(\myal+\mybt)\end{array}\right)\label{eq:Ls-notation}\end{equation}
and define a linear operator $\mK=\Id+\myde^{-1}\mLs$, where $\Id$
is the identity map in~$\RR^{2}$.
\begin{lem}
\label{lem:V-d-rec} The following recurrent equations for $\emV n$
and $\emd n$ hold\[
\emd n=\ko\left(\emd{n-1}+\myb\dt_{n}\right)\]
\[
\emV n=\mK\left(\emV{n-1}+\dt_{n}^{2}\mqd+\dt_{n}\mqo+\dt_{n}\emd{n-1}\mqz\right)\]
where $\dt_{n}=\tau_{n}-\tau_{n-1}$~.
\end{lem}
\textbf{Proof of Lemma~\ref{lem:V-d-rec}.} The dynamics of the process
$x(t)$ consists of two parts: free motion and pairwise interaction
between nodes of the network. Namely, the\emph{ interaction} is possible
only at random time moments $0<\tau_{1}<\tau_{2}<\cdots$ and has
the form of \emph{synchronizing jumps}: at time $\tau_{n}$ a pair
of nodes \[
\mbox{\ensuremath{(i,j)}, \ensuremath{\quad i=1,2,\ldots,N}, \ensuremath{\quad j=2,\ldots,N}, \ensuremath{\quad i\not=j},}\]
 is randomly chosen, and the value of clock~$j$ jumps to the value
of clock~$i$: \begin{equation}
\left(x_{i},x_{j}\right)\rightarrow\left(x_{i},x_{i}\right).\label{eq:pair-inter}\end{equation}
Inside intervals $(\tau_{k},\tau_{k+1})$ components of the process
$x(t)$ move according to the free evolutions~(\ref{MAD:eq:x-vt-sB})
and~(\ref{MAD:eq:x-rt}) driven by independent Brownian motions.
Recall that the pair $(i,j)$ mentioned above corresponds to a time-stamped
message sent by clock~$i$ to clock~$j$ at time~$\tau_{n}$. It
is easy to see that any pair $(1,j)$ is chosen with probability $\myal\myde^{-1}$
and any pair $(i,j)$, $i>1$, is chosen with probability $\mybt\myde^{-1}$.
Keeping in mind~(\ref{eq:pair-inter}) we introduce a family of maps
$\Mp_{(i,j)}:\,\RR^{N+1}\rightarrow\RR^{N+1},$ \[
\Mp_{(i,j)}:\,\,\left(x_{1},\ldots,x_{N+1}\right)\mapsto\left(x'_{1},\ldots,x'_{N+1}\right),\]
where $x'_{j}=x_{i},\quad x'_{k}=x_{k},\quad k\not=j,$ and define
a map-valued random variable $\rMp$ such that \begin{equation}
\sP\left\{ \rMp=\Mp_{(i,j)}\right\} =\begin{cases}
\myal\myde^{-1}, & i=1,\\
\mybt\myde^{-1}, & i>1.\end{cases}\label{eq:S=00003DS-ij}\end{equation}
So Lemma~\ref{lem:V-d-rec} will immediatelly follow from the next
two lemmas.
\begin{lem}
\label{lem:V-Sx}For any $x\in\RR^{N+1}$ \[
\sE V(\rMp x)=\mK V(x),\qquad\sE d(\rMp x)=\mko d(x).\]

\end{lem}
Denote an auxiliary stochastic process $z=z(t)$ evolving according
to the free dynamics~(\ref{MAD:eq:x-vt-sB})--(\ref{MAD:eq:x-rt}). 
\begin{lem}
\label{lem:V-z-free}For $s<t$

\[
\sE\left(d(z(t))\,|\, z(s)\right)=d(z(s))+\myb\dt\]
\[
\sE\left(V(z(t))\,|\, z(s)\right)=V(z(s))+\dt^{2}\cdot\mqd+\dt\cdot\mqo+\dt\cdot d(z(s))\mqz\]
where $\dt=t-s.$
\end{lem}
The proofs of Lemmas~\ref{lem:V-Sx} and~\ref{lem:V-z-free} are
quite straightforward and very similar to proofs of corresponding
lemmas in~\cite{manita-TVP-large-ident-eng,Man-F-T,Man-Rome}. They
are omitted here.

It should be noted that, in fact, Lemmas~\ref{lem:V-Sx} and~\ref{lem:V-z-free}
are valid under\emph{ weaker} conditions than the \emph{Markovian
assumption}. Namely, these lemma are true also for the following \textbf{\emph{semi-Markovian
assumption}}: messages are sent at epochs $\left\{ \tau_{n}\right\} _{n=0}^{\infty}$
where intervals $\left\{ \tau_{m}-\tau_{m-1}\right\} _{m=0}^{\infty}$
are independent random variables with identical continuous distribution;
at time moment $\tau_{n}$ the pair $(i,j)$=(sender,receiver) is
chosen independently according to the map-valued random variable $\rMp$
defined in~(\ref{eq:S=00003DS-ij}).

Briefly speaking, under Markovian assumption the sequence $\left\{ \tau_{n}\right\} _{n=0}^{\infty}$
is a Poisson flow but under semi-Markovian assumption $\left\{ \tau_{n}\right\} _{n=0}^{\infty}$
is a general renewal process~\cite{COX,Gned-Koval-1968}.

\subsection{Decomposition into sum over diagrams}

Having equations of Lemma~\ref{lem:V-d-rec} we can use them to express
$\emV n$ via $\emV{n-2}$ and so on. Proceeding recursively we get
a sum with large number of summands. We want to organize these summands
in some proper way to be able to evaluate $\emV n$ as an explicit
function of $\emV 0$. To do this we use a general approach known
as decomposition into sum over diagrams. We start from description
of a set of diagrams corresponding to our specific task.

\paragraph*{Admissible diagrams.}

Fix some $n$. Let us define a \emph{set of admissible diagrams} $\mathcal{G}(n)$
of order $n$. We say that an oriented graph $G$ belongs to the set
$\mathcal{G}(n)$ iff $G$ is a path \begin{equation}
G=\left(\vD_{0},\,\vD_{1},\,\ldots,\,\vD_{r}\right),\quad\quad1\leq r\leq n+1,\label{eq:G-path}\end{equation}
 with vertices $\vD_{0},\,\vD_{1},\ldots,\vD_{r}\in\mathcal{M}_{n}$
and oriented edges $\vD_{i-1}\rightarrow\vD_{i}$, $i=1,\ldots,r$,
satisfying to the following conditions.

a) Vertices of $G$ are labeled by pairs $(c,d)$ belonging to \[
\mathcal{M}_{n}=\left\{ \,\vD=(c,d)\,:\,(c,d)\in S_{C}^{(n)}\times S_{D}\cup\left\{ \left(\labt,2\right)\right\} \,\right\} \]
where $c\in S_{C}^{(n)}:=\left\{ 0,1,\ldots,n-1,n\right\} $ and $d\in S_{D}:=\left\{ 0,1,2\right\} $.

b) The initial vertex $\vD_{0}$ of any path $G$ is $\left(\labt,2\right)$.
For $i\geq1$ the vertex $\vD_{i}$ has the form $(n+1-i,\, l_{i})$,
$l_{i}\in S_{D}$. The final vertex may be \[
\vD_{r}=\begin{cases}
(n+1-r,\,0), & \quad\mbox{if \ensuremath{r\leq n}},\\
(0,d), & \quad\mbox{if }r=n+1\qquad\mbox{where }d\in S_{D}.\end{cases}\]

c) The first edge $\vD_{0}\rightarrow\vD_{1}$ of $G$ has the form
$\left(\labt,2\right)\rightarrow(n,l')$, $l'\in S_{D}$, and other
edges $\left(\vD_{i}\rightarrow\vD_{i+1}\right)$, $1\leq i\leq r-1$,
have the form $(n+1-i,\, l_{i})\rightarrow(n-i,\, l_{i+1})$ where
$l_{i+1}\leq l_{i}$. There are \emph{no edges} of the form $(m,0)\rightarrow(m-1,0)$.

For the path~(\ref{eq:G-path}) we use notation \[
\VG=\left\{ \vD_{0},\,\vD_{1},\,\ldots,\,\vD_{r}\right\} ,\qquad\EG=\left\{ \vD_{i-1}\rightarrow\vD_{i}\,|\, i=1,\ldots,r\right\} \]
for the set of vertices and the set of edges of~$G$. For reasons
which will be clear below, sometimes we will remap elements of the
set $S_{D}$ as follows \begin{equation}
\mbox{0\ensuremath{\longleftrightarrow\lo}, 1\ensuremath{\longleftrightarrow\ld}, 2\ensuremath{\longleftrightarrow\lV}.}\label{eq:agr-vert-remap}\end{equation}

\paragraph*{Contribution of a path.}

Let $\vec{\tau}=(\tau_{1},\ldots,\tau_{n})\in\mathbb{R}_{+}^{n}$
and $t>0$ are such that $0<\tau_{1}<\cdots<\tau_{n}<t$. We define
$\JJ_{G}=J_{G}(\vec{\tau},t)\in\mathbb{R}^{2}$, a \emph{contribution}
of the path~$G$ (see~(\ref{eq:G-path})), as follows \begin{eqnarray*}
J_{G}(\vec{\tau},t) & = & \rho_{\vD_{0}}\pi_{(\vD_{0},\vD_{1})}\rho_{\vD_{1}}\pi_{(\vD_{1},\vD_{2})}\cdots\pi_{(\vD_{r-1},\vD_{r})}\rho_{\vD_{r}}\,=\,\\
 & = & \left(\prod_{a\in\VG}\rho_{a}\prod_{b\in\EG}\pi_{b}\right)_{\mbox{ordered along the path}},\end{eqnarray*}
where functions $\rho_{a}=\rho_{a}(\vec{\tau},t)$ and $\pi_{b}=\pi_{b}(\vec{\tau},t)$
are \emph{contributions} of vertices and edges defined below. Namely,
we put\[
\begin{array}{rl}
 & \rho_{\vD_{0}}=1\\
 & \rho_{v}=\begin{cases}
\emV 0\,, & \quad\mbox{if \ensuremath{v=(0,2)}},\\
\emd 0\,, & \quad\mbox{if \ensuremath{v=(0,1)}},\\
1\,, & \quad\mbox{if \ensuremath{v=(0,0)}},\end{cases}\end{array}\quad\quad\begin{array}{r}
\rho_{v}=\begin{cases}
\matK\,, & \quad\mbox{if \ensuremath{v=(m,2)}},\\
\ko\,, & \quad\mbox{if \ensuremath{v=(m,1)}},\\
1\,, & \quad\mbox{if \ensuremath{v=(m,0)}},\end{cases}\\
\quad m=1,\ldots,n\,.\end{array}\]
and \[
\pi_{b}=f_{l,l'}(\tau_{m}-\tau_{m-1}),\quad m=1,\ldots n+1,\qquad\tau_{n+1}\equiv t\]
for the edge $b=(m,l)\rightarrow(m-1,l')$. The functions $f_{l,l'}$
are defined as follows: 

\begin{align*}
f_{22}(s) & =Id, & f_{20}(s) & =s^{2}\qdv+s\qod, & f_{11}(s) & =1,\\
 &  & f_{21}(s) & =s\mqz, & f_{10}(s) & =s\bB.\end{align*}

\begin{prop}
\label{pro:expan-diag}The following expansion holds

\[
\sE\left(V(x(t)\,|\,\,\left\{ \tau_{j}\right\} _{j=1}^{\infty}\right)=\sum_{G\in\mathcal{G}(\Pi_{t})}\JJ_{G}\,,\]
where $\mathcal{G}(n)$ is the set of admissible diagrams of order
$n$ (defined above) and $\JJ_{G}\in\mathbb{R}^{2}$ is the contribution
of a diagram~$G$. 
\end{prop}
The proof of Proposition~\ref{pro:expan-diag} is just a careful
development of the recurrent equations of Lemma~\ref{lem:V-d-rec}.
The only thing one should pay attention is the last interval $(\tau_{\Pi_{t}},t]$
of the total time segment~$[0,t]$. There is no synchronization jump
at time~$t$. Hence we need to apply for this interval only Lemma~\ref{lem:V-z-free}.
This explains assignment~$\rho_{\vD_{0}}=1$ for the root vertex~$\vD_{0}=(\labt,2)$.

\subsection{Evaluation of the functions}

Below we use agreement~(\ref{eq:agr-vert-remap}). For any path $G\in\mathcal{G}(n)$
denote by $n_{1}=n_{1}(G)$ the number of vertices of the form $(m,\ld)$,
$m\leq n$, and by $n_{2}=n_{2}(G)$ the number of vertices of the
form $(m,\lV)$, $m\leq n$. In is clear that $n_{1},n_{2}\geq0$
and $n_{1}+n_{2}\leq n+1$. Now we decompose the set of admissible
diagrams $\mathcal{G}(n)$ into the following nonintersecting subsets\[
\mathcal{G}(n)=\mathcal{G}_{0}(n)\cup\mathcal{G}_{2}(n)\cup\mathcal{G}_{10}(n)\cup\mathcal{G}_{11}(n),\]
\[
\mathcal{G}_{0}(n)=\left\{ G\in\mathcal{G}(n):\, n_{1}=0,\, n_{2}=n+1\right\} =\left\{ G_{0}\right\} \]
\[
\mathcal{G}_{2}(n)=\left\{ G\in\mathcal{G}(n):\, n_{1}=0,\,0\leq n_{2}\leq n\right\} \]
\[
\mathcal{G}_{10}(n)=\left\{ G\in\mathcal{G}(n):\,1\leq n_{1}\leq n+1,\, n_{2}=n+1-n_{1}\right\} \]
\[
\mathcal{G}_{11}(n)=\left\{ G\in\mathcal{G}(n):\,1\leq n_{1}\leq n,\, n_{2}\geq0,\, n_{1}+n_{2}<n+1\right\} \]
The subset $\mathcal{G}_{0}(n)$ consists of a single path \[
G_{0}=\left(\left(\labt,\lV\right),\,\left(n,\lV\right),\,\left(n-1,\lV\right),\ldots,\,\left(0,\lV\right)\right)\]
 and contribution of this path is \[
J_{G_{0}}(\vec{\tau},t)=\prod_{b\in\EG}\pi_{b}\prod_{a\in\VG}\rho_{a}\,\left|\begin{array}{c}
\,\\
\,\\
\!\!\,_{G=G_{0}}\end{array}\right.=\matK^{n}\emV 0.\]
For $G\in\mathcal{G}_{2}(n)$ the final vertex $\vD_{r}$ is $\left(n-n_{2},\lo\right)$
and \[
J_{G}(\vec{\tau},t)=\mK^{n_{2}}\dt_{n+1-n_{2}}^{2}\mqd+\mK^{n_{2}}\dt_{n+1-n_{2}}\mqo\,,\]
where $\dt_{k}=\tau_{k}-\tau_{k-1}$. If $G\in\mathcal{G}_{10}(n)$
then $\vD_{r}=\left(0,\ld\right)$, $r=n+1$, and\[
J_{G}(\vec{\tau},t)=\mK^{n+1-n_{1}}\mko^{n_{1}-1}\dt_{n_{1}}\emd 0\mqz\]

For $G\in\mathcal{G}_{11}(n)$ the final vertex $\vD_{r}$ is $\left(n-n_{1}-n_{2},\lo\right)$
and\[
J_{G}(\vec{\tau},t)=\mK^{n_{2}}\mko^{n_{1}}\dt_{n+1-n_{2}}\dt_{n+1-n_{1}-n_{2}}\myb\mqz\,\]

So we have \[
\sE\left(V(x(t)\,|\,\,\left\{ \tau_{j}\right\} _{j=1}^{\infty}\right)=\left(\sum_{G\in\mathcal{G}_{0}(\Pi_{t})}+\sum_{G\in\mathcal{G}_{2}(\Pi_{t})}+\sum_{G\in\mathcal{G}_{10}(\Pi_{t})}+\sum_{G\in\mathcal{G}_{11}(\Pi_{t})}\right)\JJ_{G}\,,\]
taking conditional expectation $\sE\left(\cdot\,|\,\Pi_{t}\right)$
of the both sides of this decomposition we obtain \begin{eqnarray*}
\sE\left(V(x(t)\,|\,\,\Pi_{t}=n\,\right) & = & \left(\sum_{G\in\mathcal{G}_{0}(n)}+\sum_{G\in\mathcal{G}_{2}(m)}+\sum_{G\in\mathcal{G}_{10}(n)}+\sum_{G\in\mathcal{G}_{11}(n)}\right)\JJ_{G}\,=\,\\
 & = & \matK^{n}\emV 0+\sum_{n_{2}=0}^{n}\left(\mK^{n_{2}}\Ys 2{n+1-n_{2}}\mqd+\mK^{n_{2}}\Ys 1{n+1-n_{2}}\mqo\right)\,+\\
 &  & \,+\left(\sum_{n_{2}=0}^{n}\mK^{n_{2}}\mko^{n-n_{2}}\Ys 1{n+1-n_{2}}\right)\,\emd 0\mqz+\,\\
 &  & \,+\left(\sum_{n_{2}=0}^{n-1}\,\sum_{n_{1}=1}^{n-n_{2}}\,\mK^{n_{2}}\mko^{n_{1}}\Ys{1,1}{n+1-n_{2},\, n+1-n_{1}-n_{2}}\right)\myb\mqz\end{eqnarray*}

where \begin{eqnarray*}
\Ys mk & = & \Ys mk(t,n)=\sE\left(\left(\dt_{k}\right)^{m}\,|\,\Pi_{t}=n\right),\qquad m=1,2,\\
\Ys{1,1}{i,j} & = & \Ys{1,1}{i,j}(t,n)=\sE\left(\dt_{i}\dt_{j}\,|\,\Pi_{t}=n\right),\quad i\not=j\,,\\
\dt_{j} & = & \tau_{j}-\tau_{j-1}\,.\end{eqnarray*}
 It is easy to check that in general \emph{semi-Markov case }(Subsection~\ref{sub:many-notation})
\begin{eqnarray*}
\Ys m{k_{1}}(t,n) & = & \Ys m{k_{2}}(t,n),\qquad1\leq k_{1},k_{2}\leq n,\\
\Ys{1,1}{i_{1},j_{1}}(t,n) & = & \Ys{1,1}{i_{2},j_{2}}(t,n),\quad\quad i_{1},j_{1},i_{2},j_{2}\in\left\{ 1,\ldots,n\right\} ,\quad i_{1}\not=j_{1},\, i_{2}\not=j_{2}.\end{eqnarray*}
Under \emph{Markovian} \emph{assumption} we have much stronger result:\begin{eqnarray*}
\Ys 1k(t,n) & = & \frac{t}{n+1}\,,\qquad1\leq k\leq n+1,\\
\Ys 2{k_{1}}(t,n) & = & \frac{2t^{2}}{\left(n+1\right)\left(n+2\right)}\,,\qquad1\leq k\leq n+1,\\
\Ys{1,1}{i,j}(t,n) & = & \frac{t^{2}}{\left(n+1\right)\left(n+2\right)}\,,\quad\quad i,j\in\left\{ 1,\ldots,n+1\right\} ,\quad i\not=j\,.\end{eqnarray*}

Note that $2\mqd=\myb\mqz$ and \[
\mathcal{G}_{2}(n)\cup\mathcal{G}_{11}(n)=\left\{ G\in\mathcal{G}(n):\, n_{1}\geq0,\, n_{2}\geq0,\, n_{1}+n_{2}\leq n\right\} .\]
Hence we just proved the following statement. 
\begin{lem}
In the Markovian case \textup{\begin{eqnarray}
\sE\left(V(x(t)\,|\,\,\Pi_{t}=n\,\right) & = & \matK^{n}\emV 0+\frac{t}{n+1}\left(\sum_{n_{2}=0}^{n}\mK^{n_{2}}\right)\mqo\,+\label{eq:dec-V-Kn-1}\\
 &  & \,+\frac{t}{n+1}\left(\sum_{n_{2}=0}^{n}\mK^{n_{2}}\mko^{n-n_{2}}\right)\,\emd 0\mqz+\,\label{eq:dec-V-Kn-2}\\
 &  & \,+\frac{t^{2}}{\left(n+1\right)\left(n+2\right)}\left(\sum_{n_{1}+n_{2}\leq n}\,\mK^{n_{2}}\mko^{n_{1}}\right)\myb\mqz\,\,.\label{eq:dec-V-Kn-3}\end{eqnarray}
}
\end{lem}
Now we are going to evaluate $\sE V(x(t)$ by averaging on~$\Pi_{t}$:
\begin{eqnarray}
\sE V(x(t) & = & \sE\left(\sE\left(V(x(t)\,|\,\,\Pi_{t}\,\right)\right)=\nonumber \\
 & = & \sum_{n=0}^{+\infty}\,\frac{\left(\myde t\right)^{n}}{n!}e^{-\myde t}\,\sE\left(V(x(t)\,|\,\,\Pi_{t}=n\,\right).\label{eq:average-on-Pi}\end{eqnarray}
To do this we need two technical lemmas.
\begin{lem}
\label{lem:ident-A}For any finite-dimensional matrix~$A$ the following
identities hold \begin{eqnarray}
\sE A^{\Pi_{t}} & = & e^{-\myde t(\Id-A)}\label{eq:iden-A-1}\\
\sE\frac{\, t\, A^{\Pi_{t}+1}}{\Pi_{t}+1} & = & \frac{1}{\myde}\,\left(e^{-\myde t(\Id-A)}-e^{-\myde t}\Id\right)\label{eq:iden-A-2}\\
\sE\frac{t^{2}\, A^{\Pi_{t}+2}}{(\Pi_{t}+1)(\Pi_{t}+2)} & = & \frac{1}{\myde^{2}}\,\left(e^{-\myde t(\Id-A)}-\frac{\Id+\myde tA}{e^{\myde t}}\right)\label{eq:iden-A-3}\end{eqnarray}
 where $\left(\Pi_{t},\, t\geq0\right)$ is the Poissonian process
of intensity $\myde$.
\end{lem}
Proof of Lemma~\ref{lem:ident-A} is very easy. We omit it.
\begin{lem}
\label{lem:U-a-a}Let $\ao,\at\in\left(0,1\right)$. Consider \begin{eqnarray*}
U_{1}(\ao,\at) & := & \sE\frac{t}{\Pi_{t}+1}\sum_{n_{2}=0}^{\Pi_{t}}\ao^{n_{2}}\at^{\Pi_{t}-n_{2}},\\
U_{2}(\ao,\at) & := & \sE\frac{t^{2}}{\left(\Pi_{t}+1\right)\left(\Pi_{t}+2\right)}\sum_{n_{1}+n_{2}\leq\Pi_{t}}\,\ao^{n_{2}}\at^{n_{1}},\end{eqnarray*}
where $\left(\Pi_{t},\, t\geq0\right)$ is a Poissonian process of
intensity~$\myde>0$. Then 

for $\ao\not=\at$\begin{eqnarray*}
U_{1}(\ao,\at) & = & \myde^{-1}\left(\ao-\at\right)^{-1}\left(e^{-\left(1-\ao\right)\myde t}-e^{-\left(1-\at\right)\myde t}\right)\\
\myde^{2}\cdot U_{2}(\ao,\at) & = & \frac{1}{(1-\ao)(1-\at)}-\frac{{\displaystyle \frac{e^{-\left(1-\ao\right)\myde t}}{1-\ao}}-{\displaystyle \frac{e^{-\left(1-\at\right)\myde t}}{1-\at}}}{\ao-\at}\,\end{eqnarray*}

and for $\ao=\at=\aa$ \begin{eqnarray*}
U_{1}(\aa,\aa) & = & te^{-\left(1-\aa\right)\myde t}\,,\\
\myde^{2}\cdot U_{2}(\aa,\aa) & = & \frac{1}{(1-\aa)^{2}}-\left(\frac{1}{(1-\aa)^{2}}+\frac{\myde t}{1-\aa}\right)e^{-\left(1-\aa\right)\myde t}\,.\end{eqnarray*}

\end{lem}
Proof of Lemma~\ref{lem:U-a-a} will be given in Section~\ref{sec:app}.

Now we proceed with evaluation of (\ref{eq:average-on-Pi}). Putting~$A=\mK$
in~(\ref{eq:iden-A-1}) and applying averaging~(\ref{eq:average-on-Pi})
to the first summand in~(\ref{eq:dec-V-Kn-1}) we get \[
\sE\left(\matK^{\Pi_{t}}\emV 0\right)\,=e^{-\myde t(\Id-\mK)}\sE\emV 0\,=\, e^{t\mLs}\sE\emV 0,\]
since $\mK=\Id+\myde^{-1}\mLs$ (see notation~(\ref{eq:notation-scalars})--(\ref{eq:Ls-notation})
at the beginning of Subsection~\ref{sub:many-notation}). Similarly,
using~(\ref{eq:iden-A-2}) for $A=\mK$ and $A=\Id$~, we obtain\begin{eqnarray*}
\sE\,\frac{t}{\Pi_{t}+1}\left(\sum_{n_{2}=0}^{\Pi_{t}}\mK^{n_{2}}\right)\mqo & = & \sE\,\frac{t}{\Pi_{t}+1}(\Id-\mK)^{-1}\left(\Id-\mK^{\Pi_{t}+1}\right)\mqo\\
 & = & \left(-\mLs\right)^{-1}\left(\Id-e^{\mLs t}\right)\,\mqo.\end{eqnarray*}

To find expectation of summands~(\ref{eq:dec-V-Kn-2}) and~(\ref{eq:dec-V-Kn-3})
we shall analyze spectrum of the operator $\mLs:\,\RR^{2}\rightarrow\RR^{2}$
and then apply Lemma~\ref{lem:U-a-a}. It is easy to check that operator
$\mLs$ has two different eigenvalues $\lon$ and $\ldn$, corresponding
to eigenvectors $\eon$ and $\edn$, \begin{equation}
\lon=-\myal,\qquad\ldn=-2\left({\displaystyle \myal+\mybt}\right),\label{eq:lon-ldn-def}\end{equation}
\begin{equation}
\eon=\left(\begin{array}{c}
1\\
\,\\
{\displaystyle \frac{2\myal}{{\displaystyle \myal+2\mybt}}}\end{array}\right)=\left(\begin{array}{c}
1\\
\,\\
{\displaystyle \frac{2\lon}{\ldn-\lon}}\end{array}\right),\qquad\quad\edn=\left(\begin{array}{c}
0\\
1\end{array}\right)\,.\label{eq:e1-e2-def}\end{equation}
where $\myal$ and $\mybt$ are the same as in~(\ref{eq:notation-scalars}).
Hence actions of the operators $\mko\Id$ and $\mK$ of the vector
$\eon$ are identical,\[
\mko\Id\,\eon=\mK\eon=\aa\eon,\quad\quad\aa=1-\myde^{-1}\myal=1+\myde^{-1}\lon,\]
but their actions on the vector~$\edn$ are different:\[
\mko\Id\,\edn=\ao\edn,\quad\quad\ao=1-\myde^{-1}\myal=1+\myde^{-1}\lon\,,\]
\[
\mK\edn=\at\edn,\qquad\at=1-2\myde^{-1}\left(\myal+\mybt\right)=1+\myde^{-1}\ldn\,,\qquad\ao\not=\at.\]
In the basis $\eon,\,\edn$ we have $\onull=\eon+\wdn\edn$~, where
\begin{equation}
\wdn=-{\displaystyle {\displaystyle \frac{2\myal}{{\displaystyle \myal+2\mybt}}}}\,=\,-{\displaystyle \frac{2\lon}{\ldn-\lon}}\,,\label{eq:w2-n-def}\end{equation}
Substituting this decomposition into~(\ref{eq:dec-V-Kn-2})--(\ref{eq:dec-V-Kn-3}),
we apply Lemma~\ref{lem:U-a-a} to calculate expectations of (\ref{eq:dec-V-Kn-2})
and (\ref{eq:dec-V-Kn-3}) separately on each linear subspace $\left\langle \eon\right\rangle $
and $\left\langle \edn\right\rangle $. Finally, noting that $\left(1-\aa_{i}\right)\myde=-\lin$
we come to the following statement.
\begin{prop}
\label{pro:R-D-expl-form}The functions $\MADR_{N}(t)$ and $\MADD_{N}(t)$
can be given in the following explicit form:

\begin{eqnarray*}
\left(\begin{array}{c}
\MADR_{N}(t)\\
\MADD_{N}(t)\end{array}\right) & = & e^{\mLs t}\left(\begin{array}{c}
\MADR_{N}(0)\\
\MADD_{N}(0)\end{array}\right)\,\,+\left(-\mLs\right)^{-1}\left(\Id-e^{\mLs t}\right)\,\left(\begin{array}{c}
\sigma^{2}\\
2\sigma^{2}\end{array}\right)+\,\\
 &  & \,+(v-r)\,\mbd(0)\left(te^{\lon t}\,\eon\,+\,\frac{e^{\lambda_{2,N}t}-e^{\lon t}}{\lambda_{2,N}-\lon}\,\wdn\,\edn\right)+\,\\
 &  & \,+\left(\frac{1}{\lon^{2}}-\left(\frac{1}{\lon^{2}}-\frac{t}{\lon}\right)e^{\lon t}\right)\,(v-r)^{2}\,\eon+\,\\
 &  & \,+\left(\frac{1}{\lon\lambda_{2,N}}+\frac{{\displaystyle \frac{e^{\lambda_{2,N}t}}{\lambda_{2,N}}}-{\displaystyle \frac{e^{\lon t}}{\lon}}}{\lambda_{2,N}-\lon}\right)\,\wdn\,(v-r)^{2}\,\edn\,,\end{eqnarray*}
where $\mbd(t):=\sE d(x(t))$ with $d(\cdot)$ defined in~(\ref{eq:d-x-def}). \end{prop}
\begin{rem}
\label{rem:param-N-sim} Note that many parameters in the above formula
depend on~$N$. But in the limit $N\rightarrow\infty$ these dependencies
become rather simple: \[
\mLs\sim\frac{1}{N}\,\mM=\,\frac{1}{N}\left(\begin{array}{cc}
-\alpha & \,\,0\\
2\alpha & \,\,-2(\alpha+\beta)\end{array}\right),\quad\eon\sim\left(\begin{array}{c}
1\\
{\textstyle {\displaystyle \frac{2\alpha}{\alpha+2\beta}}}\end{array}\right),\quad\edn=\left(\begin{array}{c}
0\\
1\end{array}\right),\]
\[
\lon=-\frac{\alpha}{N}\,,\quad\ldn\sim-\frac{2(\alpha+\beta)}{N}\,,\quad\wdn\sim-\frac{2\alpha}{\alpha+2\beta}\,.\]

\end{rem}
~\global\long\def\leqRD{\lefteqn{\left(\begin{array}{c}
\MADR_{N}(t)\\
\MADD_{N}(t)\end{array}\right)=e^{\mLs t}\left(\begin{array}{c}
\MADR_{N}(0)\\
\MADD_{N}(0)\end{array}\right)\,\,+\left(-\mLs\right)^{-1}\left(\Id-e^{\mLs t}\right)\,\left(\begin{array}{c}
\sigma^{2}\\
2\sigma^{2}\end{array}\right)+\,}}

\begin{rem}
\label{rem:R-D-g1-g2} Proposition~\ref{pro:R-D-expl-form} can be
rewritten in the following way: \begin{eqnarray*}
\leqRD\\
 &  & \,+(v-r)\, t\,\mbd(0)\left(\go'\left(\lon t\right)\,\eon\,+\,\frac{\go\left(\lambda_{2,N}t\right)-\go\left(\lon t\right)}{(\lambda_{2,N}-\lon)\, t}\,\wdn\,\edn\right)+\,\\
 &  & \,+(v-r)^{2}\, t^{2}\,\Bigl(\gd'\left(\lon t\right)\,\eon\,+\,\frac{\gd\left(\lambda_{2,N}t\right)-\gd\left(\lon t\right)}{(\lambda_{2,N}-\lon)\, t}\,\wdn\,\edn\,\Bigr),\end{eqnarray*}
where $\go(y)=e^{y},\quad\gd(y)={\displaystyle \frac{e^{y}-1}{y}}\,.$
Using~(\ref{eq:w2-n-def}) we easily transform this representation
as follows \begin{eqnarray*}
\leqRD\\
 &  & \,+(v-r)\, t\,\mbd(0)\left(\go'\left(\lon t\right)\,\onull\,+\,\Fn{\go}t\,\zo\right)+\,\\
 &  & \,+(v-r)^{2}\, t^{2}\,\Bigl(\gd'\left(\lon t\right)\,\onull\,+\,\Fn{\gd}t\,\zo\Bigr),\end{eqnarray*}
where \[
\Fn{\gG}t=\left(\frac{\gG\left(\lambda_{2,N}t\right)-\gG\left(\lon t\right)}{(\lambda_{2,N}-\lon)\, t}\,-\gG'\left(\lon t\right)\right)\,\wdn\,\,.\]

\end{rem}
For the function $\mbd(t)$ we also have explicit formula.
\begin{prop}
\label{pro:d-N-expl-form}The function $d_{N}(t)$ has the following
form \begin{equation}
\mbd(t)=\mbd(0)\exp\left(-\frac{\alpha}{N}\, t\right)+\left(1-\exp\left(-\frac{\alpha}{N}\, t\right)\right)\,\frac{(v-r)N}{\alpha}\,\label{eq:dN-t-expl}\end{equation}
and hence is a solution to the following equation \begin{equation}
\frac{d}{dt}\,\mbd(t)=-\frac{\alpha}{N}\,\mbd(t)+(v-r).\label{eq:dN-t-equation}\end{equation}

\end{prop}
Proof of this proposition can be obtained by the same method as the
proof of Proposition~\ref{pro:R-D-expl-form} but the corresponding
reasonings are much more shorter and simpler. So we omit the proof
of Proposition~\ref{pro:d-N-expl-form} and simply refer to~\cite{Man-F-T}
where similar statement was presented in full details.

The following statement follows from Proposition~\ref{pro:R-D-expl-form}
by direct calculation.
\begin{cor}
$(\MADR_{N}(t),\MADD_{N}(t))$ is a solution of the following system
of differential equations:

\begin{eqnarray*}
\frac{d}{dt}\left(\begin{array}{c}
\MADR_{N}(t)\\
\MADD_{N}(t)\end{array}\right) & = & \left(\begin{array}{cc}
-\alpha/N & 0\\
2\alpha/N & -2(\alpha/N+\beta/(N-1)\,)\end{array}\right)\left(\begin{array}{c}
\MADR_{N}(t)\\
\MADD_{N}(t)\end{array}\right)+\\
 &  & \,+\left(\begin{array}{c}
2(v-r)\mbd(t)\\
0\end{array}\right)+\left(\begin{array}{c}
\sigma^{2}\\
2\sigma^{2}\end{array}\right)\end{eqnarray*}
where the function $d_{N}(t)$ is the same as in~(\ref{eq:dN-t-expl})
and~(\ref{eq:dN-t-equation}).
\end{cor}

\subsection{Study of the asymptotical behavior}

In this subsection we obtain all results on asymptotic behavior which
were presented in Section~\ref{sec:main-res-WSN-serv}.

\subsubsection*{Proof of Theorem~\ref{thm:N-fixed-t-to-inf}}

Here $N$ is fixed and $t$ tends to infinity. Since matrix~$\mLs$
has two different negative eigenvalues one easily concludes that in
explicit representation of Proposition~\ref{pro:R-D-expl-form} all
terms containing $e^{\mLs t}$ or $e^{\lin t}$ go to zero as $t\rightarrow\infty$.
So the limit is equal to\[
\left(\begin{array}{c}
\MADR_{N}(\infty)\\
\MADD_{N}(\infty)\end{array}\right)=\left(-\mLs\right)^{-1}\left(\begin{array}{c}
\sigma^{2}\\
2\sigma^{2}\end{array}\right)+(v-r)^{2}\,\left(\frac{1}{\lon^{2}}\,\eon+\frac{1}{\lon\lambda_{2,N}}\,\wdn\,\edn\,\right).\]
Matrix $\left(-\mLs\right)^{-1}$ can be calculated explicitly:\[
\left(-\mLs\right)^{-1}=\left(\begin{array}{cc}
\,\,\myal^{-1}\,\, & 0\\
(\myal+\mybt)^{-1}\, & \,\frac{1}{2}(\myal+\mybt)^{-1}\,\end{array}\right).\]
Using (\ref{eq:e1-e2-def}) and~(\ref{eq:w2-n-def}) after some algebraic
transformation we get\[
\MADR_{N}(t)\rightarrow\MADR_{N}(\infty)=\sigma^{2}\myal^{-1}+\frac{\left(v-r\right)^{2}}{\lon^{2}}\,=\frac{\sigma^{2}}{\alpha}\, N+\frac{\left(v-r\right)^{2}}{\alpha^{2}}\, N^{2}\,,\]
\[
\MADD_{N}(t)\rightarrow\MADD_{N}(\infty)=\frac{2\,\sigma^{2}}{\myal+\mybt}\,+\frac{2\,\left(v-r\right)^{2}}{\lon\ldn}\,\sim\,\,\frac{2\,\sigma^{2}}{\alpha+\beta}\, N\,+\frac{\left(v-r\right)^{2}}{\alpha\cdot\left(\alpha+\beta\right)}\, N^{2}\,.\]
This proves Theorem~\ref{thm:N-fixed-t-to-inf}.

\subsubsection*{Proofs of Theorems~\ref{thm:v-eq-r-3-phases} and~\ref{thm:v-NOT-eq-r-3-phases}.}

At the beginning we study the first time scale, namely, we assume
that \begin{equation}
\tn\rightarrow\infty\,\mbox{ but }\tn/N\rightarrow0\,\mbox{ as }N\rightarrow\infty\,.\label{eq:1st-time-scale}\end{equation}
 To do this we shall use the next lemma.
\begin{lem}
\label{lem:g-yN-zN}Consider a function $\gG=\gG(y)$ which assumed
to be analytical in some neighborhood of the point $y=0$. Assume
also that $\gG''(0)\not=0$. Then for any sequences $\left\{ y_{N}\right\} $
and $\left\{ z_{N}\right\} $ tending to~$0$ as $N\rightarrow\infty$
we have \[
\gG'\left(y_{N}\right)\,\eon\,+\,\frac{\gG\left(z_{N}\right)-\gG\left(y_{N}\right)}{z_{N}-y_{N}}\,\wdn\,\edn=\,\]
\[
\,=\gG'\left(y_{N}\right)\oz+\,\left(\frac{1}{2}\,\gG''(0)\left(z_{N}-y_{N}\right)\wdn+\left(z_{N}-y_{N}\right)O\left(\left|z_{N}\right|+\left|y_{N}\right|\right)\right)\,\left(\begin{array}{c}
0\\
1\end{array}\right)\,,\]
where the notation $a_{N}=O(b_{N})$ means, as usual, that $\,{\displaystyle \limsup_{N}\left|\frac{a_{N}}{b_{N}}\right|}<+\infty.$
\end{lem}
Proof of this lemma is a straightforward calculation. 

Since the both $\lon$ and $\ldn$ are negative and condition~(\ref{eq:R0D0})
holds we remark that the first summand in the r.h.s.~of Proposition~\ref{pro:R-D-expl-form}
is bounded. Under assumption~(\ref{eq:1st-time-scale}) \[
\lon\tn\rightarrow0\,\mbox{ and\,\,}\ldn\tn\rightarrow0\,.\]
So the second summand has a rather simple behavior:\begin{equation}
\left(-\mLs\right)^{-1}\left(\Id-e^{\mLs\tn}\right)\,\left(\begin{array}{c}
\sigma^{2}\\
2\sigma^{2}\end{array}\right)\sim\left(\begin{array}{c}
\sigma^{2}\\
2\sigma^{2}\end{array}\right)\,\tn+\left(\begin{array}{c}
O(\tn/N)\\
O(\tn/N)\end{array}\right)\,\tn\,.\label{eq:1st-tsc-1}\end{equation}
The next summands of the representation in Proposition~\ref{pro:R-D-expl-form}
demands very careful analysis. We use Remark~\ref{rem:R-D-g1-g2}
and Lemma~\ref{lem:g-yN-zN} for $y_{N}=\lon\tn$ and $z_{N}=\ldn\tn$.
\[
\go'\left(\lon t\right)\,\eon\,+\,\frac{\go\left(\lambda_{2,N}t\right)-\go\left(\lon t\right)}{(\lambda_{2,N}-\lon)\, t}\,\wdn\,\edn=\]
\[
=e^{y_{N}}\oz+\,\left(\frac{1}{2}\,\left(z_{N}-y_{N}\right)\wdn+\left(z_{N}-y_{N}\right)O\left(\tn/N\right)\right)\,\left(\begin{array}{c}
0\\
1\end{array}\right)\,=\]
\[
=\left(1+O(\lon\tn)\,\right)\oz+\,\left(\,\mlon\tn+O\left(\tn^{2}/N^{2}\right)\right)\,\left(\begin{array}{c}
0\\
1\end{array}\right)\,=\]
\begin{equation}
=\left(\begin{array}{c}
1+O(\tn/N)\,\\
\mlon\tn+O\left(\tn^{2}/N^{2}\right)\end{array}\right)\label{eq:1st-tsc-2}\end{equation}
where we have use identity $\left(z_{N}-y_{N}\right)\wdn=2\mlon\tn$
(see~(\ref{eq:lon-ldn-def}) and~(\ref{eq:w2-n-def})). Similarly,
\[
\gd'\left(\lon t\right)\,\eon\,+\,\frac{\gd\left(\lambda_{2,N}t\right)-\gd\left(\lon t\right)}{(\lambda_{2,N}-\lon)\, t}\,\wdn\,\edn=\]
\[
=\frac{1-(1-y_{N})e^{y_{N}}}{y_{N}^{2}}\,\oz+\,\left(\frac{1}{2}\cdot\frac{1}{3}\left(z_{N}-y_{N}\right)\wdn+\left(z_{N}-y_{N}\right)O\left(\tn/N\right)\right)\,\left(\begin{array}{c}
0\\
1\end{array}\right)\,=\]
\[
=\left(\frac{1}{2}+O(\lon\tn)\,\right)\oz+\,\left(\,\frac{1}{3}\mlon\tn+O\left(\tn^{2}/N^{2}\right)\right)\,\left(\begin{array}{c}
0\\
1\end{array}\right)\,=\]
\begin{equation}
=\left(\begin{array}{c}
\frac{1}{2}+O(\tn/N)\,\\
\frac{1}{3}\mlon\tn+O\left(\tn^{2}/N^{2}\right)\end{array}\right)\label{eq:1st-tsc-3}\end{equation}
Putting (\ref{eq:1st-tsc-1})--(\ref{eq:1st-tsc-3}) into formula
of Remark~\ref{rem:R-D-g1-g2} we get \begin{eqnarray*}
\MADR_{N}(\tn) & = & \sigma^{2}\tn+O(\tn/N)+\,(v-r)\,\mbd(0)\,\tn+(v-r)\tn O(\tn/N)+\\
 &  & \,+\frac{1}{2}(v-r)^{2}\tn^{2}+(v-r)^{2}\tn^{2}O(\tn/N)\end{eqnarray*}
\begin{eqnarray*}
\MADD_{N}(\tn) & = & 2\sigma^{2}\tn+O(\tn/N)+\,\\
 &  & \,+(v-r)\,\mbd(0)\,\mlon\tn^{2}+(v-r)\tn O\left(\tn^{2}/N^{2}\right)+\,\\
 &  & \,+\frac{1}{3}(v-r)^{2}\,\mlon\tn^{3}+(v-r)^{2}\tn^{2}O\left(\tn^{2}/N^{2}\right)\end{eqnarray*}
In the case $v-r=0$ these formulae immediately imply item~P1 of
the Theorem~\ref{thm:v-eq-r-3-phases}. 

Consider the case $v-r\not=0$ . Since $\tn\rightarrow\infty$ terms
containing $(v-r)\tn$ are asymptotically smaller than corresponding
terms containing $(v-r)^{2}\tn^{2}$. Moreover, $\mlon\tn^{2}=\tn O(\tn/N)$.
Hence \begin{eqnarray*}
\MADR_{N}(\tn) & = & \sigma^{2}\tn+O(\tn/N)+\,(v-r)\,\mbd(0)\,\tn+\,\\
 &  & \,+\frac{1}{2}(v-r)^{2}\tn^{2}+(v-r)^{2}\tn^{2}O(\tn/N)\end{eqnarray*}
\begin{eqnarray}
\MADD_{N}(\tn) & = & 2\sigma^{2}\tn+O(\tn/N)+(v-r)\,\mbd(0)\,\tn O(\tn/N)+\,\nonumber \\
 &  & \,+\frac{1}{3}(v-r)^{2}\,\mlon\tn^{3}+(v-r)^{2}\tn^{2}O\left(\tn^{2}/N^{2}\right).\label{eq:DN-1st-sc-tmp}\end{eqnarray}
So our next conclusion is that under assumptions~(\ref{eq:1st-time-scale})
and $v-r\not=0$ \[
\MADR_{N}(\tn)\sim\frac{1}{2}(v-r)^{2}\tn^{2}\,.\]
To analyze $\MADD_{N}(\tn)$ we should compare $\tn$ and $\mlon\tn^{3}$.
Recalling that $\mlon=\myal=\alpha/N$ we split the time scale~(\ref{eq:1st-time-scale})
into three subscales\[
\frac{\tn}{\sqrt{N}}\,\rightarrow0,\qquad\frac{\tn}{\sqrt{N}}\,\rightarrow\co,\quad\co>0,\qquad\mbox{and}\quad\frac{\tn}{\sqrt{N}}\,\rightarrow\infty\quad\mbox{but}\quad\frac{\tn}{N}\,\rightarrow0.\]
Considering~(\ref{eq:DN-1st-sc-tmp}) on each subscale we get that
\begin{itemize}
\item if $\frac{\tn}{\sqrt{N}}\,\rightarrow0$, then $\MADD_{N}(\tn)\sim2\sigma^{2}\tn$,
\item if $\frac{\tn}{\sqrt{N}}\,\rightarrow\co,\quad\co>0,$ then \begin{eqnarray*}
\MADD_{N}(\tn) & \sim & 2\sigma^{2}\tn+\frac{1}{3}\alpha(v-r)^{2}\tn^{3}/N\,\sim\,\\
 & \sim & 2\sigma^{2}\co\sqrt{N}+\frac{1}{3}\alpha(v-r)^{2}\,\co^{3}\sqrt{N}\sim\left(2\sigma^{2}+\frac{1}{3}\alpha(v-r)^{2}\,\co^{2}\right)\tn\end{eqnarray*}

\item if $\frac{\tn}{\sqrt{N}}\,\rightarrow\infty$ but $\frac{\tn}{N}\,\rightarrow0$,
then $\MADD_{N}(\tn)\sim\frac{1}{3}\alpha(v-r)^{2}\,\tn^{3}/N$.
\end{itemize}
Next let us prove Theorems~\ref{thm:v-eq-r-3-phases} and~\ref{thm:v-NOT-eq-r-3-phases}
for the time scale \begin{equation}
\tn/N\rightarrow c>0,\quad\quad N\rightarrow\infty.\label{eq:tN-N-c}\end{equation}
Taking into account Remark~\ref{rem:param-N-sim} and assumption~(\ref{eq:tN-N-c})
we get\[
\tn\mLs\rightarrow c\mM,\quad\quad\lon\tn\rightarrow-\alpha c,\quad\quad\ldn\tn\rightarrow-2(\alpha+\beta)c\]
where $\mM$ is defined in~(\ref{eq:matr-M}). If $\gG$ is an analytic
function then $\gG'(\lon\tn)=\gG'(-\alpha c)+o(1)$ and $\Fn{\gG}{\tn}=\hH(\gG,c)+o(1),$
where \[
\hH(\gG,c):=\left(\frac{\gG\left(-2(\alpha+\beta)c\right)-\gG\left(-\alpha c\right)}{-2(\alpha+\beta)c+\alpha c}\,-\gG'\left(-\alpha c\right)\right)\cdot\left(-\frac{2\alpha}{\alpha+2\beta}\right)\,.\]

\global\long\def\vnutrRDcM{\left(\begin{array}{c}
 \MADR_{N}(\tn)\\
\MADD_{N}(\tn)\end{array}\right)\sim(-cM)^{-1}\left(\Id-e^{cM}\right)\,\left(\begin{array}{c}
 \sigma^{2}\\
2\sigma^{2}\end{array}\right)\,\tn\,+o(\tn)\,+\,}

\global\long\def\leqRDcM{\lefteqn{\vnutrRDcM}}

Using Remark~\ref{rem:R-D-g1-g2} we have \begin{eqnarray*}
\leqRDcM\\
 &  & \,+(v-r)\,\tn\,\sE\mbd(0)\left(\go'(-\alpha c)\,\onull\,+\,\hH(\go,c)\,\zo+\left(\begin{array}{c}
o(1)\\
o(1)\end{array}\right)\right)+\,\\
 &  & \,+(v-r)^{2}\,\tn^{2}\,\Bigl(\left(\gd'(-\alpha c)+o(1)\right)\,\onull\,+\,\left(\hH(\gd,c)+o(1)\right)\,\zo\Bigr).\end{eqnarray*}
Now if $v=r$ we easily get item~P2 of Theorem~\ref{thm:v-eq-r-3-phases}.
In the case $v\not=r$ we obtain item~P2 of Theorem~\ref{thm:v-NOT-eq-r-3-phases}
with\begin{eqnarray*}
\cdho(c) & = & \gd'(-\alpha c),\\
\cdhd(c) & = & \hH(\gd,c),\end{eqnarray*}
where $\gd$ is the same as in Remark~\ref{rem:R-D-g1-g2}. Now we
are able to get explicit forms of the functions $\cdho(c)$ and $\cdhd(c)$:
\begin{eqnarray}
\cdho(c) & = & \frac{1-(1+\alpha c)e^{-\alpha c}}{\alpha^{2}c^{2}}\,,\nonumber \\
\cdhd(c) & = & 2\, c^{-2}\cdot\frac{1-(1+\alpha c)e^{-\alpha c}}{\alpha\left(\alpha+2\beta\right)}-\,\label{eq:h1-h2-expl-form}\\
 &  & \,-{\displaystyle {\displaystyle \frac{2\, c^{-2}\alpha}{{\displaystyle \alpha+2\beta}}}\left(\frac{1}{2\left(\alpha+\beta\right)\alpha}-\frac{1}{\alpha+2\beta}\left({\displaystyle \frac{e^{-\alpha c}}{\alpha}}-{\displaystyle \frac{e^{-2(\alpha+\beta)c}}{2(\alpha+\beta)}}\right)\right)}\,.\nonumber \end{eqnarray}

The study of the time scale $\tn/N\rightarrow\infty$ is very similar
to the proof of Theorem~\ref{thm:N-fixed-t-to-inf} (see above).
So we omit details.\label{sec-proofs-Ths}

\subsection{Time scale analysis}

\label{sub:time-sc-analysis}

In this subsection we derive some \emph{corollaries} from our main
theorems. We examine behavior of the functions $\MADR_{N}(\tn)$ and
$\MADD_{N}(\tn)$ on special time scales $\tn=\mys N^{\gamma}$ where
$\myg>0$ and $\mys>0$. Summing up Theorems~\ref{thm:N-fixed-t-to-inf}--\ref{thm:v-NOT-eq-r-3-phases}
and Subsection~\ref{sec-proofs-Ths} we get the following general
statement: Assume that $N\rightarrow\infty$ and assumption~(\ref{eq:R0D0})
holds. Then \begin{eqnarray*}
\MADR_{N}(\mys N^{\gamma}) & \sim & \cR(\mys,\gamma)\, N^{\psR(\gamma)}\,,\\
\MADD_{N}(\mys N^{\gamma}) & \sim & \cD(\mys,\gamma)\, N^{\psD(\gamma)}\,,\end{eqnarray*}
where all functions $\cR,\cD,\psR,\psD$ are positive. Moreover, these
functions can be calculated explicitely: 
\begin{description}
\item [{Case~1:}] $v=r$ --- zero skew.

\begin{tabular}{ccl}
$\psR(\myg)=\min(\myg,1)$, & $\qquad$ & $\cR(\mys,\myg)=\begin{cases}
\sigma^{2}\mys, & \gamma<1\,,\\
\sigma^{2}\mylR(\mys), & \gamma=1,\\
\sigma^{2}/\alpha,\quad & \gamma>1,\end{cases}$\tabularnewline
$\psD(\myg)=\min(\myg,1),$ &  & $\cD(\mys,\myg)=\begin{cases}
2\sigma^{2}\mys, & \gamma<1\,,\\
2\sigma^{2}\mylD(\mys), & \gamma=1,\\
2\sigma^{2}/(\alpha+\beta),\quad & \gamma>1,\end{cases}$\tabularnewline
\end{tabular}

where \begin{eqnarray*}
\mylR(\mys) & = & \gd(-\alpha\mys)\,,\\
\mylD(\mys) & = & \frac{\alpha\,\gd(-\alpha\mys)+2\beta\,\gd(-2(\alpha+\beta)\mys)}{\alpha+2\beta}\,\end{eqnarray*}
and the function $\gd$ is the same as~Remark~\ref{rem:R-D-g1-g2}.

\item [{Case~2:}] $v\not=r$ --- nonzero skew. 
\end{description}
\noindent \begin{tabular}{ccl}
$\psR(\myg)=\min(2\myg,2)$, & \null ~~ & $\cR(\mys,\myg)=\begin{cases}
\frac{1}{2}\left(v-r\right)^{2}\mys^{2}, & \gamma<1\,,\\
\left(v-r\right)^{2}\mys^{2}\cdho(\mys), & \gamma=1,\\
{\displaystyle \frac{\left(v-r\right)^{2}}{\alpha^{2}}}\,,\quad & \gamma>1,\end{cases}$\tabularnewline
 &  & \tabularnewline
$\psD(\myg)=\begin{cases}
\gamma, & \gamma\leq\frac{1}{2},\\
3\gamma-1, & \lefteqn{{\textstyle \frac{1}{2}}<\gamma<1,}\\
2, & \gamma\geq1,\end{cases}$ &  & $\cD(\mys,\myg)=\begin{cases}
2\sigma^{2}\mys, & \gamma<\frac{1}{2},\\
2\sigma^{2}\mys+\frac{1}{3}\alpha(v-r)^{2}\,\mys^{3}, & \gamma=\frac{1}{2}\\
\frac{1}{3}\alpha(v-r)^{2}\,\mys^{3}, & \frac{1}{2}<\gamma<1,\\
(v-r)^{2}\,\mys^{2}\cdhd(\mys), & \gamma=1,\\
(v-r)^{2}\alpha^{-1}\left(\alpha+\beta\right)^{-1}, & \gamma>1.\end{cases}$\tabularnewline
\end{tabular}\\
where $\cdho$ and $\cdhd$ are the same as in~(\ref{eq:h1-h2-expl-form}).

~

How to interpret and to explain the above results? The choice of the
time scale $\tn=\mys N^{\gamma}$ means that we consider a new\emph{
time unit} which is equal to $N^{\gamma}$ units of the physical time~$t$.
Then paramemer~$s$ is the time measured in new units. We see that
behavior of the stochastic system $x(t)=\left(x_{1}(t),\ldots,x_{N+1}(t)\right)$
for large $N$ is different for different values of~$\myg$. We recall
(Section~\ref{sec:st-model-Ns}) that $\MADR_{N}(\mys N^{\gamma})$
and $\MADD_{N}(\mys N^{\gamma})$ are squares of clock synchronization
errors. The above corollaries show existence of several phases in
the evolution of the network. In the case $v=r$ there are three such
phases, they were discussed in Remark~\ref{rem:v=00003Dr-stages}. 

In the case $v\not=r$ on {}``small times'' $(\myg<1)$ the function
$\MADR_{N}(\mys N^{\gamma})$ is proportional to the square of~$\tn$.
We can imagine that on these times effect of clock synchronization
is negligible and the asymptotical value $\MADR_{N}(\mys N^{\gamma})$
depends only on summands of the form \begin{equation}
\left(x_{1}(t)-x_{j}(t)\right)^{2}=\left(x_{1}(0)-x_{j}(0)+(r-v)t\right)^{2}\sim(v-r)^{2}t^{2}.\label{eq:x1-xj-t2}\end{equation}
 For the scale $\myg=1$ the total number of synchronization messages
from the server~1 is of order~$N$, their influence become important
and we observe a slowdown: for $\mys\rightarrow\infty$\[
\mys^{2}\cdho(\mys)\,\uparrow\,\alpha^{-2}\,<\,+\infty.\]
On the scales $\myg>1$ we see the result of full strength competition
between synchronization jumps and the desynchronization generated
by motions of the client nodes $2,\ldots,N+1$. This is the final
synchronization phase. 

In in the same case $v\not=r$ the behavior of the function~$\MADD_{N}(\mys N^{\gamma})$
is much more interesting. It appears that the {}``small times''
$(\myg<1)$ are splitted into three subphases: $0<\myg<\frac{1}{2}$~,
$\myg=\frac{1}{2}$ and $\frac{1}{2}<\myg<1$. Recall that the function~$\MADD_{N}$
describes the internal inconsistency of the client's clocks $x_{2},\ldots,x_{N+1}$.
On the first subphase ($0<\myg<\frac{1}{2}$) influence of the server
node~1 is negligible in comparison with a {}``noise'' produced
by identical client nodes, hence we see the same asymptotics $2\sigma^{2}\mys\, N^{\myg}$
as in the Case~1 (see also asymptotics in~\cite{Man-F-T}). On the
scale $\myg=\frac{1}{2}$ these two forces are of the same order $N^{1/2}$
and $\cD(\mys,\frac{1}{2})=2\sigma^{2}\mys+\frac{1}{3}\alpha(v-r)^{2}\mys^{3}$.
On the third subphase ($\frac{1}{2}<\myg<1$) the server node~1 dominates
over the free dynamics of the clients. Surprisingly, that for scales
$\frac{1}{2}<\myg<1$ the influence of the server of the accurate
time produces a desynchronization (but not synchronization) of the
client's clocks. The explanations is the following one: on these scales
only \emph{small part} of client nodes had interacted with the node~1
till the time~$\mys N^{\gamma}$. The most part of client nodes (of
order~$N$) still have no idea about the server's clock~$x_{1}$.
Adjusting of a client's clock to the value $x_{1}$ brings a \emph{large
number} of summands of the form~(\ref{eq:x1-xj-t2}) to the function~$\MADD_{N}$.
On the scale $\myg=1$ the number of nodes interacted with the server~1
is of order~$N$ and again we observe a slowdown of desynchronization.
The time scales $\gamma>1$ correspond to the final synchronization
phase.

\section{Appendix}

\label{sec:app}

\paragraph*{Proof of Lemma~\ref{lem:U-a-a}. }

Let $\ao\not=\at$. It is straightforward to check that \[
\sum_{n_{2}=0}^{n}\ao^{n_{2}}\at^{n-n_{2}}=\left(\ao-\at\right)^{-1}\left(\ao^{n+1}-\at^{n+1}\right).\]
Using~(\ref{eq:iden-A-2}) for $A=\ao$ and $A=\at$ we find $U_{1}(\ao,\at)$.

Now let us calculate $U_{2}(\ao,\at)$. After some simple algebra
we have identity\[
\sum_{\begin{array}{c}
n_{1}\ge0,\, n_{2}\geq0\\
n_{1}+n_{2}\leq n\end{array}}\,\ao^{n_{1}}\at^{n_{2}}=\frac{{\displaystyle \frac{1}{1-\ao}}-{\displaystyle \frac{1}{1-\at}}}{\ao-\at}\,-\,\frac{{\displaystyle \frac{\ao^{n+2}}{1-\ao}}-{\displaystyle \frac{\at^{n+2}}{1-\at}}}{\ao-\at}\,,\]
Applying~(\ref{eq:iden-A-3}) for $A=1$, $A=\ao$ and $A=\at$~,
we get that $U_{2}(\ao,\at)$ \emph{multiplied by} $\myde^{2}$ is
equal to \[
\frac{{\displaystyle \frac{1}{1-\ao}}-{\displaystyle \frac{1}{1-\at}}}{\ao-\at}\,\left(1-\frac{1+\myde t}{e^{\myde t}}\right)\,-\,\]
\[
\,-\,\frac{{\displaystyle \frac{\left(e^{-\myde t(1-\ao)}-{\displaystyle \frac{1+\myde t\ao}{e^{\myde t}}}\right)}{1-\ao}}-{\displaystyle \frac{\left(e^{-\myde t(1-\at)}-{\displaystyle \frac{1+\myde t\at}{e^{\myde t}}}\right)}{1-\at}}}{\ao-\at}\]
By direct transformations and cancellation of terms this form can
be reduced to the following one\[
\,-\frac{{\displaystyle \frac{e^{-\myde t(1-\ao)}-1}{1-\ao}}-{\displaystyle \frac{e^{-\myde t(1-\at)}-1}{1-\at}}}{\ao-\at}=\frac{1}{(1-\ao)(1-\at)}-\frac{{\displaystyle \frac{e^{-\myde t(1-\ao)}}{1-\ao}}-{\displaystyle \frac{e^{-\myde t(1-\at)}}{1-\at}}}{\ao-\at}\,\,.\]
This proves the statement of Lemma~\ref{lem:U-a-a} for $\ao\not=\at$. 

The case $\ao=\at$ is simpler than just considered case $\ao\not=\at$.
So we omit details here. Note, that explicit expressions for $U_{i}(\aa,\aa)$,
$i=1,2$, correspond to \emph{formal limits} of $U_{i}(\ao,\at)$
as $\ao\rightarrow\aa$, $\at\rightarrow\aa$. $\Box$

\section{Conclusion and future work}

\label{sec:future-w}

We proposed a basic probabilistic model of clock synchronization in
large WSNs interacting with an accurate time server. It was shown
that in large networks ($N\rightarrow\infty)$ there exists several
time scales~$t=\tn$ of qualitatively different collective behavior
of the network. In other words, a large network passes different phases
on its road to synchronization. The \emph{phase of effective syn\-chro\-ni\-zation}
is the most interesting among them. For our basic network we give
a detailed description of this phase. Moreover, explicit formulae
obtained in~Theorems~\ref{thm:N-fixed-t-to-inf}--\ref{thm:v-NOT-eq-r-3-phases}
provide keys to future analytical study of various optimization and
performance evaluation problems related to WSNs.

We believe that obtained results about existence of several different
phases in evolution should take place for more general classes of
\emph{large} networks. Mathematical tools used in the present study
will work also for general nonhomogeneous time synchronization model.
Of course, in that case we cannot expect such short and explicit results
as in Section~\ref{sec:main-res-WSN-serv}. Successful examples of
synchronization studies for some special nonhomogeneous (or weakly
nonhomogeneous) systems can be found in~\cite{man-Shch,mal-man-TVP,Francois-LNCS,manita-umn}.

Our basic model can be developped also for more general and realistic
assumptions about message sending algorithms. We plan to get rid of
condition~(\ref{eq:exp-distrub}) to be able to consider arbitrary
distributed intervals between messages. The theory of general random
flows~\cite{COX} can be useful here but, unfortunately, the corresponding
probabilistic model will be non-Markovian (similarly to~\cite{Man-Rome}).
The white noise assumption in the clock model~(\ref{MAD:eq:x-vt-sB})
is confusing for the following reason: it contradicts to a general
clock modelling principle that time never run backward (see page~313
in~\cite{SundBuyKshem-review}). To justify our choice we note, first,
that from the mathematical viewpoint this assumption is convenient
but not necessary and it will be removed in future papers, and, in
the second place, this is an usual assumption for many modern clock
models~(see~\cite{Li-Rus,Zhu-Ma-Ryu-pat2009,Aus4auth2011})

We hope that presented approach combined with other methods will be
useful for analysis of many practical distributed systems. 

\end{document}